\documentclass[12pt,leqno]{article}
\usepackage{amsfonts}
\usepackage{amsmath, amsthm, amsfonts, amssymb, color}
\usepackage{mathrsfs}
\usepackage{color}
\theoremstyle{definition}

\newcommand{\scr}[1]{\mathscr #1}
\definecolor{wco}{rgb}{0.5,0.2,0.3}

\numberwithin{equation}{section} \theoremstyle{remark}

\newcommand{\ua}{\uparrow}

\title{{\bf   Wasserstein Convergence Rate   for  Empirical Measures of   Markov Processes } }
\author{{\bf  Feng-Yu Wang   }\\
\footnotesize{ Center for Applied Mathematics, Tianjin University, Tianjin 300072, China}\\
\footnotesize{    wangfy@tju.edu.cn}}
\begin{document}
\allowdisplaybreaks
\def\R{\mathbb R}  \def\ff{\frac} \def\ss{\sqrt} \def\B{\mathbf
B}
\def\N{\mathbb N} \def\kk{\kappa} \def\m{{\bf m}}
\def\ee{\varepsilon}\def\ddd{D^*}
\def\dd{\delta} \def\DD{\Delta} \def\vv{\varepsilon} \def\rr{\rho}
\def\<{\langle} \def\>{\rangle}
  \def\nn{\nabla} \def\pp{\partial} \def\E{\mathbb E}
\def\d{\text{\rm{d}}} \def\bb{\beta} \def\aa{\alpha} \def\D{\scr D}
  \def\si{\sigma} \def\ess{\text{\rm{ess}}}\def\s{{\bf s}}
\def\beg{\begin} \def\beq{\begin{equation}}  \def\F{\scr F}
\def\Ric{\mathcal Ric} \def\Hess{\text{\rm{Hess}}}
\def\e{\text{\rm{e}}} \def\ua{\underline a} \def\OO{\Omega}  \def\oo{\omega}
 \def\tt{\tilde}\def\[{\lfloor} \def\]{\rfloor}
\def\cut{\text{\rm{cut}}} \def\P{\mathbb P} \def\ifn{I_n(f^{\bigotimes n})}
\def\C{\scr C}      \def\aaa{\mathbf{r}}     \def\r{r}
\def\gap{\text{\rm{gap}}} \def\prr{\pi_{{\bf m},\nu}}  \def\r{\mathbf r}
\def\Z{\mathbb Z} \def\vrr{\nu} \def\ll{\lambda}
\def\L{\scr L}\def\Tt{\tt} \def\TT{\tt}\def\II{\mathbb I}
\def\i{{\rm in}}\def\Sect{{\rm Sect}}  \def\H{\mathbb H}
\def\M{\mathbb M}\def\Q{\mathbb Q} \def\texto{\text{o}} \def\LL{\Lambda}
\def\Rank{{\rm Rank}} \def\B{\scr B} \def\i{{\rm i}} \def\HR{\hat{\R}^d}
\def\to{\rightarrow} \def\gg{\gamma}
\def\EE{\scr E} \def\W{\mathbb W}
\def\A{\scr A} \def\Lip{{\rm Lip}}\def\S{\mathbb S}
\def\BB{\scr B}\def\Ent{{\rm Ent}} \def\i{{\rm i}}\def\itparallel{{\it\parallel}}
\def\g{{\mathbf g}}\def\Sect{{\mathcal Sec}}\def\T{\mathcal T}\def\BB{{\bf B}}
\def\f\ell \def\g{\mathbf g}\def\BL{{\bf L}}  \def\BG{{\mathbb G}}
\def\Bd{{D^E}} \def\BdP{D^E_\phi} \def\Bdd{{\bf \dd}} \def\Bs{{\bf s}} \def\GA{\scr A}
\def\Bg{{\bf g}}  \def\Bdd{\psi_B} \def\supp{{\rm supp}}\def\div{{\rm div}}
\def\ddiv{{\rm div}}\def\osc{{\bf osc}}\def\1{{\bf 1}}\def\BD{\mathbb D}
\def\H{{\bf H}}\def\gg{\gamma} \def\n{{\mathbf n}}\def\GG{\Gamma}\def\HAT{\hat}
\def\SU{{\bf SU}}\def\Var{{\rm Var}}\def\DD{\Delta}\def\T{\mathbb T}
\maketitle

\begin{abstract}  The convergence rate in Wasserstein distance  is estimated  for    empirical measures   of  ergodic Markov processes, and  the estimate can be sharp in some specific situations.
   The main result is applied to  subordinations of  typical models excluded by existing results, which include:    stochastic Hamiltonian systems on $\R^{n}\times \R^{m}$,      spherical velocity Langevin processes on  $\R^n\times\mathbb S^{n-1},$  multi-dimensional Wright-Fisher type diffusion processes, and    stable type jump processes.

  \end{abstract} \noindent

 \vskip 2cm

 \section{Introduction}

 The purpose of this paper is to provide  a general result  on  the Wasserstein convergence rate of empirical measures, which applies to a broad class of ergodic Markov  processes including    typical   models beyond the range of existing results.

\subsection{Problem in existing study}
 In recent years the Wasserstein convergence rate has been intensively investigated for the empirical measures of continuous time  stochastic systems, see  \cite{W1, W3,W4,W2,MT} for  symmetric diffusion processes,    \cite{LW1, LW2,LW3,LW4,W23,WW} for   subordinate diffusion processes,   \cite{HT,LW4} for  the fractional Brownian motion on flat torus. See also \cite{Zhu} for the study of weighted empirical measures
of symmetric diffusions on compact manifolds.

In  these references, the symmetric part of the generator has discrete spectrum with positive spectral gap, see  \cite[(2.6)]{W23}. In particular,   there exists a constant $c>0$ such that the following Poincar\'e inequality
holds:
\beq\label{FM} \mu(f^2)-\mu(f)^2 \le c\, \hat \EE(f,f),\ \ \hat\EE(f,f):=  - \mu(fLf).\end{equation}
However, this restriction excludes   degenerate models  where $\EE(f,f)$ is reducible.

For instance,   consider  the following stochastic Hamiltonian system
 $X_t:=(X_t^{(1)}, X_t^{(2)})$ on $\R^{d}\times \R^d=\R^{2d}$:
  \beq\label{E0} \beg{cases} \d X_t^{(1)}= X_t^{(2)}\d t,\\
\d X_t^{(2)}= \ss 2 \,\d W_t -\big\{\nn V(X_t^{(1)})+X_t^{(2)}\big\}\d t,\end{cases}\end{equation}
where $W_t$ is the Brownian motion on $\R^d$,  and $V\in C^2(\R^d)$   such that $Z_V:=\int_{\R^d}\e^{-V(x)}\d x<\infty$ and
$ \|\nn^2 V\|\le C(1+|\nn V|)$   holds or some constant $C>0$.
In this case,  the solution of \eqref{E0} is a diffusion process having invariant probability measure
 $\mu=\mu_V\times\scr N_1,$  where
$\mu_V(\d x)= Z_V^{-1} \e^{-V(x)}\d x$  and   $\scr N_1$ is the standard Gaussian measure on $\R^d$.
The energy form associated with \eqref{E0} is
$$\hat \EE(f,f)  =\mu\big(\big|\nn_{x^{(2)}} f(x^{(1)}, \cdot)\big|^2(x^{(2)})\big),\ \ f\in C_0^2(\R^{2d}).$$
Since the gradient is only taken for the second variable, when $f$ is a non-constant function only depending on $x^{(1)}$, the Poincar\'e inequality \eqref{FM} does not hold.

On the other hand, the Markov process $X_t=(X_t^{(1)}, X_t^{(2)})$ solving \eqref{E0} may be exponential ergodic (see \cite{V}), so it is natural to ask for the convergence rate of the empirical measure with respect to a reasonable Wasserstein distance.

 \subsection{New idea of the present work}

To derive sharp estimates on   Wasserstein distance of empirical measures, we need to  regularize the empirical measures such that analytic inequalities apply. In previous references, the empirical
 measures  are regularized  by using  the semigroup $\hat P_t$ generated by the symmetric part of the   Markov generator under study,  for which we need to make assumptions on $\hat P_t$.

For instance, let $\hat L:=\DD+\nn V$ on a compete connected Riemannian manifold $M$ for a smooth function $V$ such that $\mu(\d x):=\e^{V(x)}\d x$ is a probability measure, where $\d x$ stands for the Riemannian volume measure, and let $\W_p$ be the $p$-Wasserstein distance induced by the Riemannian distance $\rr$, see \eqref{WP} below. Then by \cite[Theorem 2]{Ledoux}, for any  probability measure $\nu(\d x):=f(x)\mu(\d x)$ on $M$,
we have
\beq\label{Ledoux} \W_p(\nu,\mu)\le p\Big(\mu\big(|\nn (-\hat L)^{-1}(f-1)|^p\big) \Big)^{\ff 1 p}.\end{equation}
  See also \cite{AMB} for a refined estimate for $\W_2$.  Since the empirical measure of the  diffusion process generated by $\hat L$ is singular with respect to $\mu$, to apply the estimate \eqref{Ledoux}, one regularizes the empirical measure  by the diffusion semigroup $\hat P_t$, see \cite{AMB} and  \cite{W1}-\cite{W2}.

  However, the above technique   does not apply to degenerate models like
 stochastic Hamiltonian systems arising from kinetic mechanics, where  the  symmetric part $\hat L$ of the generator does not induce any distance.

 To overcome this problem, we choose a different symmetric semigroup $\hat P_t$, which is not generated by the symmetric part of the underlying  Markov process, but has the same invariant probability measure $\mu.$   By choosing such a symmetric diffusion semigroup satisfying conditions needed in the study, we are able to apply \eqref{Ledoux} for the generator $\hat L$ of $\hat P_t$ to derive explicit convergence rate
 of the empirical measure with respect to the Wasserstein distance, which is  induced by the intrinsic distance of $\hat L$.

 \subsection{Organization of the paper}

 In Section 2, we introduce the framework of  the present study,   state the main result for exponential ergodic Markov processes (Theorem \ref{T}),
  and an extension for non-exponential ergodic Markov processes  (Theorem \ref{TE} ).    The convergence  rate presented in Theorem \ref{T} is  sharp for specific models  shown by    Examples 2.1-2.2 and Remark 6.1. Theorem \ref{TE} applies to any Markov process  whose semigroup converges to the invariant probability measure at certain rate corresponding  to the weak Poincar\'e inequality introduced in \cite{RW01}.

 In Sections 3-6, we apply the main result to subordinations of several  typical models:     stochastic Hamiltonian systems,  spherical velocity Langevin processes,   Wright-Fisher type diffusion processes, and   stable like processes. These models arise  from different applied areas, and are not covered by existing results on Wasserstein convergence rate of empirical measures.

\section{Framework and main result}

Let  $(M,\rr)$ be a  length space,   let $\scr P(M)$ be the set of all probability measures on $M$. For any
$p\in [1,\infty)$, the $L^p$-Wasserstein distance is defined as
\beq\label{WP}\W_p(\nu_1,\nu_2):= \inf_{\pi\in \C(\nu_1,\nu_2)} \bigg(\int_{M\times M} \rr(x,y)^p\pi(\d x,\d y)\bigg)^{\ff 1 p},
\ \ \nu_1,\nu_2\in \scr P(M),\end{equation}
where   $\C(\nu_1,\nu_2)$ is the  set of all couplings for $\nu_1$ and $\nu_2$.
We study the convergence rate in $\W_p$ for empirical measures of  ergodic Markov processes on $M$.

\subsection{Subordinate  Markov process}

 Let $X_t$ be a standard time-homogenous  Markov process on $M$  having   invariant probability measure $\mu\in \scr P(M).$  The associated Markov semigroup $P_t$ is defined as
$P_tf(x):=\E^x[f(X_t)]$ for $ t\ge 0, f\in \B_b(M),$
where $\B_b(M)$ is  the class of bounded measurable functions on $M$,      $\E^x$ is  the expectation taken
for the underlying Markov process starting at point $x$. In general, for any $\nu\in \scr P(M)$,
 $ \E^\nu $  denotes  the expectation for the  Markov process with initial distribution $\nu$.

An important class of   Markov jump processes are the subordinations (time changes) of diffusion processes  induced by Bernstein functions.
A typical model is the $\aa$-stable process generated by the fractional Laplacian, which is the time change of Brownian motion and has been used as L\'evy noise in SDEs.
See the monograph \cite{CF} and references therein for the study of subordinated Markov processes and applications.

To make time changes (i.e. the subordination)  of the Markov process $X_t$, we introduce the class {\bf B} of Bernstein functions $B$
   with $B(0)=0.$ Recall that  a Bernstein function is a function $B\in C([0,\infty))\cap C^\infty((0,\infty))$   satisfying
  $  (-1)^{n-1}\ff{\d^n B(s)}{\d s^n} \ge 0,\ \ s>0.$
  For any   $\aa\in (0,1]$,  let
\beg{align*} {\bf B}^\aa:=\Big\{B\in {\bf B}:\ \liminf_{r\to\infty}B(r)r^{-\aa}>0\Big\}.\end{align*}
Obviously, ${\bf B}^0={\bf B}.$

For each $B\in {\bf B}$, there exists a unique  stable increasing process $S_t^B$ on $[0,\infty)$ with Laplace transform
\beq\label{LT} \E\big[\e^{-r S_t^B}\big]= \e^{-B(r)t},\ \ \ t,r\ge 0,\end{equation}
see for instance \cite{SSV}.
Let $S_t^B$ be independent of $X_t$. Consider the subordinate diffusion process
$ X_t^B:= X_{S_t^B},  t\ge 0,$ and its empirical measures
$$\mu_t^B:=\ff 1 t \int_0^t \dd_{X_s^B}\d s,\ \ t>0.$$
We investigate the convergence rate of $\W_p(\mu_t^B,\mu)\to 0$ as $t\to\infty$.

\subsection{Reference symmetric diffusion process}

Let $\hat X_t$ be a reversible Markov process on $M$ with   the same   invariant probability measure $\mu,$  and with $\rr$ as the intrinsic distance. Heuristically,
$\hat X_t$ has symmetric Dirichlet form $(\hat \EE,\D(\hat\EE))$ in $L^2(\mu)$ satisfying
$$ \hat\EE(f,f)= \int_M |\nn f|^2\d\mu,\ \ f\in C_{b,L}(M)\subset \D(\hat\EE),$$
where $C_{b,L}(M)$ be the set of all bounded Lipschitz continuous functions on $M$, and
$$|\nn f(x)|:=\limsup_{y\to x} \ff{|f(y)-f(x)|}{\rr(x,y)},\ \ x\in M.$$
More precisely, we assume that $C_{b,L}(M)$ is a dense subset of $\D(\hat\EE)$
 and
$$\hat\EE(f,g)=\int_M \GG(f,g)\d\mu,\ \ f,g\in C_{b,L}(M),$$ where
$$\GG: C_{b,L}(M)\times C_{b,L}(M)\to \B_b(M)$$ is a symmetric  local square field (champ de carr\'e), i.e. for
   any $f,g,h\in C_{b,L}(M)$ and $\phi\in C_b^1(\R),$
\beg{align*} &\ss{\GG(f,f)(x)}=|\nn f(x)|:=\limsup_{y\to x} \ff{|f(y)-f(x)|}{\rr(x,y)},\ \ x\in M,\\
&\GG(fg,h)= f\GG(g,h)+ g\GG(f,h),\ \ \ \GG(\phi(f), h)=  \phi'(f) \GG(f,h).\end{align*}
Moreover,  the generator $(\hat L,\D(\hat L))$ satisfies the chain rule
$$\hat L\phi(f)= \phi'(f)\hat Lf +\phi''(f) |\nn f|^2,\ \ \ f\in \D(\hat L)\cap C_{b,L}(M), \phi\in C^2(\R).$$

\subsection{Main result}

Note that for $B(r)=r$, we have $X_t^B=X_t$ so that $\mu_t^B$ reduces to the empirical measure $\mu_t:=\ff 1 t\int_0^t  \dd_{X_s}\d s$.
In this paper, we aim to estimate
  $\W_p(\mu_t^B,\mu)$ for a   general Bernstein function $B$, which includes  $\W_p(\mu_t,\mu)$ as special case. To this end,   we   make the following assumption.

  For any $p\ge q\ge 1$, let $\|\cdot\|_{L^q(\mu)\to L^p(\mu)}$ be the operator norm from $L^q(\mu)$ to $L^p(\mu).$
  Let $(\hat P_t)_{t\ge 0}$ be the semigroup of the reversible Markov process $\hat X_t$, i.e.
  $$\hat P_t f(x):= \E\big(f(\hat X_t)\big|\hat X_0=x\big),\ \ t\ge 0, f\in \B_b(M),$$
  where $\B_b(M)$ is the set of all bounded measurable functions on $M$,

\emph{\beg{enumerate} \item[$(A_1)$] Let $p\in [2,\infty)$.  $\hat P_t$ has heat kernel $\hat p_t$ with respect to $\mu$, and there exist constants $\bb,\ll,d,k\in (0,\infty)$ such that
\beg{align}
     \label{A1}&\|\nn   \hat P_t\|_{L^2(\mu)\to L^p(\mu)}\le   k\e^{-\ll t} t^{-\bb} ,\   \ t>0,\\
    \label{A3} &\int_M  \hat P_t\rr(x,\cdot)^p (x)\mu(\d x)  \le (k t)^{\ff p 2},\ \ t\in (0,1], x\in M,\\
  \label{A4} &\int_M  \hat p_t(x,x) \mu(\d x)\le k  (1\land  t)^{-\ff d 2},\ \ t >0,\\
   \label{A5} &\|P_t-\mu\|_{L^2(\mu)}  \le k\e^{-\ll t},\ \ t\ge 0.\end{align}\end{enumerate}}
Note that in $(A_1)$   the only condition we need for the Markov process is \eqref{A5},   while other conditions \eqref{A1}-\eqref{A4} are made for the reference semigroup $\hat P_t$
which is flexible in applications. Indeed,
for smaller distance $\rr$, the energy form $\hat\EE$ is bigger, so that $\hat P_t$ has better properties. For instance, let $\mu$ be a probability measure on a connected Riemannian manifold comparable with   the volume measure,     when   the Riemmanian distance  $\rr$ is small enough we have large enough Dirichlet form $\hat\EE(f,f):=\mu(|\nn f|^2)$ such that  $ {\rm gap}(\hat L)>0,$   see \cite{CW} where   the stronger log-Sobolev inequality
is considered.

Condition \eqref{A3} refers to the $\ff 1 2$-H\"older continuity of the symmetric diffusion process $\hat X_t$, which is true for a broad class of diffusion processes. Indeed, for
a diffusion process $\hat X_t$ with generator satisfying
$$\hat L \rr(x,\cdot)^p\le c \big(1+\rr (x,\cdot)^p\big)$$
for some constant $c>0$, which is the case when $\hat X_t$ solves an SDE on $\R^d$ with linear growth coefficients and $\rr(x,y):=|x-y|$,
we have
$$\hat P_t\rr(x,\cdot)^p(x) =\E^x[\rr(x,\hat X_t)^p] \le c\int_0^t \e^{cs}\d s \le c\e^c t,\ \ t\in [0,1],$$
which implies \eqref{A3} for some constant $k>0$ and any probability measure $\mu$.
 Condition \eqref{A4} is a standard upper bound estimate on the heat kernel   for $d$-dimensional elliptic diffusions, see for instance \cite[Theorem 2.3.6]{Davies}.  Moreover,  according to Proposition \ref{P1} below,
   when $p=2$, condition \eqref{A1} with   $\bb=\ff 1 2$ follows from the existence of spectral gap, i.e.
   $${\rm gap}(\hat L):=\inf\big\{\hat\EE(f,f):\ f\in \D(\hat\EE), \mu(f)=0, \mu(f^2)=1\big\}>0.$$
   In this case, \eqref{A1} follows from \eqref{A4}, since the later  implies that $\hat L$ has discrete spectrum and hence has a spectral gap, see \cite[Theorem 3.3.19]{W05}.
When $p>2$, we have
$$\|\nn\hat P_t\|_{L^2(\mu)\to L^p(\mu)} =\|\nn\hat P_{\ff t 2} \hat P_{\ff t 2}\|_{L^2(\mu)\to L^p(\mu)} \le \|\nn \hat P_{\ff t 2}\|_{L^2(\mu)\to L^2(\mu)}\|\hat P_{\ff t 2}\|_{L^2(\mu)\to L^p(\mu)},$$
so that   \eqref{A1}  follows from ${\rm gap}(\hat L)>0$ together with a suitable upper bound of $\|\hat P_t\|_{L^2(\mu)\to L^p(\mu)},$  which is available for elliptic diffusions on compact
manifolds, see the proof of Example 2.2 for details.

    \

Before moving on, let us compare the above conditions with those in \cite[$(A_1)$]{W23}:
  there exist  constants $c,\ll >0, d\ge d'\ge 1$ and a map $k: (1,\infty)\to (0,\infty)$ such that
\beq\label{0A10} \|\hat P_t-\mu\|_{1\to\infty}\le c   t^{- \ff d 2} \e^{-\ll t},\ \ t>0,\end{equation}
\beq\label{0A1-1} \ll_i\ge c i^{\ff 2 {d'}},\ \ \ i\in \mathbb Z_+,\end{equation}
\beq\label{0A120} |\nn   \hat P_t f| \le   k(p)  ( \hat P_t |\nn f|^p)^{\ff 1 p},\ \ t\in [0,1],  p\in (1,\infty), f\in C_{b,L}(M).\end{equation}

The first essential difference is that $\hat P_t$ in  \cite[$(A_1)$]{W23} is associated with the symmetric part of the Dirichlet form for the underlying Markov process, while
$\hat P_t$ in the present framework is essentially different, the only link between the present $\hat P_t$ and the underlying Markov semigroup $P_t$ is that they share the invariant probability measure $\mu$.

Since   the generator $\hat L$ of $\hat P_t$ is not the symmetric part of the generator $L$ for  the   studied  Markov process, the eigenvalues of   $\hat L$ has nothing to do with the behavior of the underlying Markov process,
so the condition \eqref{0A1-1} is dropped from the present assumption $(A_1)$.

As we will use $\hat P_t$ to regularize the empirical measures, we adopt the conditions \eqref{A1} and  \eqref{A4} for  the gradient   and heat kernel estimates on $\hat P_t$,
where  \eqref{A1} is comparable with \eqref{0A10} for small time.
Again, because the eigenvalues of $\hat L$ has nothing to do with the underlying Markov generator, the spectral representation of
$\hat P_t$ is no longer useful for the study,   we need the gradient estimate   \eqref{A1} rather than \eqref{0A120}, where the later is easier to verify in applications.

  \beg{prp}\label{P1} If ${\rm gap}(\hat L)>0$, then for any $\ll\in (0,{\rm gap}(\hat L))$ there exists a constant $k>0$ such that
  $$\|\nn \hat P_tf\|_{L^2(\mu)}\le k   t^{-\ff 1 2} \e^{-\ll t}\|f\|_{L^2(\mu)}, \  \  t>0, f\in L^2(\mu).$$\end{prp}
  \beg{proof} Denote $\ll_1:={\rm gap}(\hat L),$ let $(E_s)_{s\ge 0}$ be the spectral family of $-\hat L$. We find a constant $k>0$ such that
  \beg{align*} &\|\nn \hat P_tf\|_{L^2(\mu)}^2=\int_{\ll_1}^\infty s\e^{-2s t}\d E_s(f)\le \e^{-2\ll t}\Big(\sup_{s\ge \ll_1} s \e^{-2(s-\ll)t}\Big) \int_{\ll_1}^\infty \d E_s(f)\\
&\le k    t^{-1} \e^{-2\ll t}\int_{\ll_1}^\infty \d E_s(f)\le  k   t^{-1} \e^{-2\ll t} \|f\|_{L^2(\mu)}^2,\ \ t>0, f\in L^2(\mu).\end{align*}
    \end{proof}

We also need the  following defined quantity $d' \in (0,\infty]$ induced by $P_t$.

\beg{defn}\label{DEF}     $d'$ is the smallest positive constant such that    the heat kernel  $p_t$ of $P_t$ with respect to $\mu$ satisfies
\beq\label{A6} \int_{M\times M} p_t(x,y)^2\mu(\d x)\mu(\d y) \le  k (1\land t)^{-\ff{d'}2},\ \ t>0. \end{equation}
If $p_t$ does not exist, or $p_t$ exists but \eqref{A6} does not hold for any $d'\in (0,\infty)$, we denote $d'=\infty.$
\end{defn}

For constants $\bb,d$ in  $(A_1)$, and $d'$ in Definition \ref{DEF}, we  denote
 $$K_{\bb,d,d',\aa}:=    \bb+\ff d 8 \bigg[1+ \Big(1-\ff{4\aa}{d'}\Big)^+ \bigg], \ \ \aa\in [0,1].$$
Moreover,  for any $t>0,$ let
\beq\label{xi}  \xi_t:= \beg{cases} t^{-1}, &\text{if}\ K_{\bb,d,d',\aa}<1,\\
t^{-1}[\log(2+t)]^2, &\text{if}\ K_{\bb,d,d',\aa}=1,\ d' \ne 4\aa,\\
t^{-1}[\log(2+t)]^3, &\text{if}\ K_{\bb,d,d',\aa}=1,\ d' =4\aa,\\
t^{-\ff{1}{2K_{\bb,d,d',\aa}-1 }}, &\text{if}\ K_{\bb,d,d',\aa}>1,\ d' \ne 4\aa,\\
t^{-\ff{1}{2K_{\bb,d,d',\aa}-1 }}\log(2+t), &\text{if}\ K_{\bb,d,d',\aa}>1,\ d' = 4\aa. \end{cases}
\end{equation}

\beg{thm}\label{T} Assume $(A_1)$ for some $p\in [2,\infty)$ and let $B\in {\bf B}^\aa$ for some $\aa\in [0,1].$ Then   there exists a constant $c>0$ such that
\beq\label{WP} \E^\mu[\W_p(\mu_t^B,\mu)^2]\le  c  \xi_t,\ \ \ t>0. \end{equation}
If the semigroup $P_t^B$ of $X_t^B$ has heat kernel  $p_t^B$ with respect to $\mu$, then for any $q\in [1,2]$ and $x\in M,$
\beq\label{WP'} \E^x [\W_p(\mu_t^B,\mu)^q]\le \ff{2^{q-1}}{t^q} \int_0^1 \E^x \big[\mu(\rr(X_s^B,\cdot)^p)^{\ff q p}\big]\d s + 2^{q-1}\|p_1^B(x,\cdot)\|_{L^{\ff 2{2-q}}(\mu)} (c \xi_{t-1})^{\ff q 2},\ \ t>1.\end{equation}

\end{thm}

\paragraph{Remark 2.1.}  (1) The reference semigroup $\hat P_t$ from $(A_1)$ will be used to regularize the empirical measure $\mu_t^B$ into
$\mu_{t,r}^B:= \mu_t^B\hat P_r$ for $r\in (0,1)$, which has density with respect to $\mu$ so that the estimate \eqref{Ledoux} applies to $\W_p(\mu_{t,r}^B,\mu),$ see the proof of Theorem \ref{T} for details.
In particular, for the stochastic Hamiltonian system \eqref{E0}, conditions \eqref{A1}-\eqref{A4} hold for $\hat P_t$ generated by
$\hat L:=\DD -\nn H$ on $\R^{2d}$, where $H(x^{(1)}, x^{(2)}):= V(x^{(1)})+ \ff 1 2 |x^{(2)}|^2$ for $ x^{(1)}, x^{(2)}\in \R^d,$ see the proof of Theorem \ref{T2} below with $m=n=d,\kk=1$ and $Q$ being the  identity matrix.

(2) By the standard Markov property,   for any $\nu\in \scr P$ and $q\in [1,2]$ with $h_\nu:=\ff{\d\nu}{\d\mu}\in L^{\ff 2{2-q}}(\mu)$, \eqref{WP} implies
\beg{align*}&\big(\E^\nu[\W_p(\mu_t^B,\mu)^q]\big)^{\ff 2 q}=\bigg(\int_M h_\nu (x)  \E^x[\W_p(\mu_t^B,\mu)^q]\mu(\d x)\bigg)^{\ff 2 q}\\
&\le  \bigg(\int_M h_\nu (x)  \big(\E^x[\W_p(\mu_t^B,\mu)^2]\big)^{\ff q 2}\mu(\d x)\bigg)^{\ff 2 q}\\
&\le \|h_\nu\|_{L^{\ff 2 {2-q}}(\mu)}^{\ff 2 q}  \int_M \E^x[\W_p(\mu_t^B,\mu)^2]\mu(\d x) \\
&=\|h_\nu\|_{L^{\ff 2 {2-q}}(\mu)}^{\ff 2 q}  \E^\mu[\W_p(\mu_t^B,\mu)^2]   \le  c \|h_\nu\|_{L^{\ff 2 {2-q}}(\mu)}^{\ff 2 q}    \xi_t,\ \ \ t>0.\end{align*}

(3) It is easy to see that $\xi_t$ is decreasing in $d'$, so estimates in  Theorem \ref{T}  remain true if   $d'$ is replaced by $\infty$, for which $K_{\bb,d,d',\aa}=K:= \bb+\ff d 4$ and $\xi_t$ reduces to
\beq\label{KK} \xi_t(K):= \beg{cases} t^{-1}, &\text{if}\ K<1,\\
t^{-1} [\log(2+t)]^2, &\text{if}\ K=1,\\
t^{-\ff 1 {2K-1}}, &\text{if}\ K>1.\end{cases}\end{equation}    Therefore,
  when $\mu$ is good enough such that the associated symmetric diffusion semigroup $\hat P_t$ satisfies conditions \eqref{A1}-\eqref{A4}, then for any Markov process
satisfying \eqref{A5} and any $B\in {\bf B},$ there exists a constant $c>0$ such that
$$  \E^\mu[\W_p(\mu_t^B,\mu)^2] \le c \xi_t(K),\ \ t>0.$$

To illustrate   Remark 2.1(2),  we present below two  examples, which provide a uniform Wasserstein convergence rate for empirical measures of Markov processes with   given invariant measure $\mu$, where the uniform rate is sharp in the second example.

\paragraph{Example 2.1.}    Let $M=\R^n$,  let   $V\in C^2(\R^{n})$ such that
$ V(x)=\psi(x) + (1+\theta |x|^2)^\tau,\  x\in\R^n, $  where  $\psi\in C_b^2(\R^{n}),$ $\theta>0,\tau\in (\ff 1 2,\infty]$ are constants.
Let $ \mu(\d x)=\mu_V(\d x):= \ff{\e^{-V(x)}\d x}{\int_{\R^n}\e^{-V(x)}\d x}.$   Then for  any Markov process on $\R^n$
satisfying \eqref{A5} and any $B\in {\bf B},$ there exists a constant $c>0$ such that  \
\beq\label{1} \E^\mu[\W_2(\mu_t^B,\mu)^2] \le c \beg{cases} t^{-1}, &\text{if}\ n=1, \tau>1,\\
t^{-1} [\log(2+t)]^2, &\text{if}\ n=1, \tau=1,\\
t^{-\ff {2\tau-1}{\tau n}}, &\text{ otherwise}. \end{cases}\end{equation}

\beg{proof}   Let $\hat L=\DD-\nn V.$  By Remark 2.1(2), it suffices to verify \eqref{A1}-\eqref{A4} for $p=2,\bb=\ff 1 2,$ and $d=\ff{2\tau n}{2\tau-1}.$ Since
$ \lim_{|x|\to\infty} \hat L|\cdot|(x)=-\infty<0,$
  \cite[Corollary 1.4]{W99} ensures  ${\rm gap}(\hat L)>0$, so that by Proposition \ref{P1},
    \eqref{A1} holds  for $p=2$ and $\bb=\ff 1 2.$

 Next, by  \cite[Theorem 2.4.4]{Wbook} and $\nn^2 V\ge -c_1I_n,$ we find a constant $c_2>0$ such that
 \beq\label{AA}\hat p_r(x,x)\le \ff {c_2} {\mu(B(x,\ss r))},\ \ x\in \R^{n}, r\in (0,1],\end{equation}
 where $B(x,r):=\{y\in \R^n: |x-y|<r\}, r>0.$
Then \eqref{A4} with $d=\ff{2\tau n}{2\tau-1}$  follows provided
\beq\label{CC}  \int_{\R^{n }}\ff {\mu(\d x)} {\mu(B(x,r))}\le c r^{-\ff{2\tau n}{2\tau-1}},\ \ r\in (0,1] \end{equation} holds for some constant $c>0.$
Below we prove this estimate.

 Since  $\psi$ is bounded, there exists a constant $C>1$ such that
\beq\label{**0} C^{-1} \e^{-(1+\theta |x|^2)^\tau}\d x\le \mu(\d x)\le C \e^{-(1+\theta |x|^2)^\tau}\d x.\end{equation}    On the other hand, by
$$\ff{|x|}2\le |x|-\ff r 4 \le |x|,\ \   r\in [0,1],\ |x|\ge 1,$$  we   find a constant  $c_3>0$ such that
\beg{align*} &\Big(1+ \theta \Big|x- \ff{rx}{4|x|}\Big|^2\Big)^\tau= (1+\theta |x|^2 )^\tau +\int_0^r \ff{\d}{\d s}  \Big(1+ \theta \big(|x|- \ff{s}{4}\big)^2\Big)^\tau\d s\\
&= (1+\theta |x|^2 )^\tau - \ff{\tau\theta}2\int_0^r \Big(1+ \theta \big(|x|- \ff{s}{4}\big)^2\Big)^{\tau-1} \Big(|x|-\ff s 4\Big)\d s\\
&\le  (1+\theta |x|^2 )^\tau- c_3 r |x|^{2\tau-1},\ \  r\in [0,1],\ |x|\ge 1.\end{align*}
Hence,  there exist  constants  $c_4, c_5>0$ such that for $|x|\ge 1$ and $r\in (0,1]$,
\beq\label{SO} \beg{split} &\mu(B(x,r))  \ge c_4 \int_{B\big(x-\ff{r x}{2|x|},\ff r 4\big)}\e^{-(1+\theta  |y|^2)^\tau}\d y\\
&\ge c_5r^n  \e^{-(1+ \theta |x- \ff{rx}{4|x|}|^2)^\tau}  \ge c_5 r^{n}\e^{-(1+\theta |x|^2)^\tau + c_3 r |x|^{2\tau-1}},\ \  r\in [0,1],\ |x|\ge 1.\end{split}\end{equation}
Combining this with \eqref{**0}, we  find  constants $c_6,c_7>0$ such that
\beg{align*}& \int_{\R^{n }}\ff {\mu(\d x)} {\mu(B(x,r))}=\int_{B(0,1)} \ff {\mu(\d x)} {\mu(B(x,r))}+\int_{B(0,1)^c} \ff {\mu(\d x)} {\mu(B(x,r))} \\
&\le C^2 r^{-n} +    c_6 r^{-n} \int_{\R^n} \e^{- c_3  r  |x|^{2\tau-1}} \d x
  =C^2 r^{-n} +  c_7 r^{-n}\int_0^\infty s^{n-1} \e^{-c_3 rs^{2\tau-1}}\d s\\
&= C^2 r^{-n} +c_7 r^{-\ff{2\tau n}{2\tau-1}}\int_0^\infty s^{n-1} \e^{-c_3 s^{2\tau-1}}\d s,\ \ r\in (0,1].\end{align*}
 This implies  \eqref{CC}  for some constant $c>0.$

 Finally,   it is easy to see that  $\nn^2 V\ge -c I_n$ and $|\nn V(x)|^2\le c(1+|x|^{4\tau})$ hold for some constant $c>0$. So, we find a constant $c_8>0$ such that
\beq\label{BO} \beg{split} &\hat L|x-\cdot|^2= 2n + 2 \<\nn V, x-\cdot\> = 2n+ 2\<\nn V(x), x-\cdot\> -2\<\nn V(x)-\nn V, x-\cdot\> \\
&\le 2n + |\nn V(x)|^2 + |x-\cdot|^2 + 2c_1|x-\cdot|^2 \le c_8 (1+|x|^{4\tau}+|x-\cdot|^2),\ \ x\in \R^{n}.\end{split} \end{equation} This implies
 \beq\label{B} \hat P_t|x-\cdot|^2(x)=\E^x|x-\hat X_t|^2 \le c_8\big(1+|x|^{4\tau}  |\big)t\e^{c_8t},\ \ x\in\R^{n}, \ t>0.\end{equation}
 Noting that  $\mu( |\cdot|^{4\tau})<\infty,$   we verify condition \eqref{A3} for $p=2$ and some constant $k>0.$
 \end{proof}

In the next example,   the upper bound \eqref{2} is sharp. Indeed,  according to \cite[Corollary 1.3(2)]{WW},
for $n\ge 3$ and  $B\in {\bf B}^\aa,$ we have
$ \inf_{x\in M} \E[\W_2(\mu_t^B,\mu)^2] \ge c t^{-\ff 2{n-2\aa}}.$
With  $\aa\to 0$ this lower bound reduces to the upper bound in \eqref{2}.

\paragraph{Example 2.2.}  Let $M$ be an $n$-dimensional compact connected Riemanian manifold possibly with a boundary $\pp M$. Let $V\in C^2(M)$ such that $\mu(\d x)=\e^{-V(x)}\d x$ is a probability measure on $M$. Then for  any Markov process on $M$
satisfying \eqref{A5} and any $B\in {\bf B},$ there exists a constant $c>0$ such that  \
\beq\label{2} \E^\mu[\W_2(\mu_t^B,\mu)^2] \le c \beg{cases} t^{-1}, &\text{if}\ n=1,\\
t^{-1} [\log(2+t)]^2, &\text{if}\ n=2,\\
t^{-\ff 2 n}, &\text{if}\ n\ge 3. \end{cases}\end{equation}
In general, for any $p\ge 2$, there exists a constant $c>0$ such that
\beq\label{2'} \E^\mu[\W_p(\mu_t^B,\mu)^2] \le c \beg{cases} t^{-1}, &\text{if}\ n=1, p\in [2,3),\\
t^{-1} [\log(2+t)]^2, &\text{if}\ n(p-1)=2,\\
t^{-\ff 2 {n(p-1)}}, &\text{otherwise}. \end{cases}\end{equation}

\beg{proof} Let $\hat X_t$ be the diffusion process generated by $\hat L=\DD-\nn V$ with reflecting boundary if it exists. Since $M$ is a compact connected Riemannian manifold, \eqref{A1}-\eqref{A4} are well known for $p=2,\bb=\ff 1 2$ and $d=n$, hence \eqref{2} holds according to Remark 2.1(2). Moreover, it is classical that
$ \|\hat P_t\|_{L^2(\mu)\to L^p(\mu)}\le c (1\land t)^{-\ff {n(p-2)}{4p}} $ holds for some constant $c>0 $ and all $t>0.$  Then
  \eqref{A1} holds for  $p\ge 2$ and $\bb=\ff 1 2+ \ff{n(p-2)}{4p}$. By Remark 2.1(2), \eqref{2'} holds.
\end{proof}

\subsection{Proof of Theorem \ref{T}}
We will apply the  estimate \eqref{Ledoux} for $p\in [2,\infty).$
 This inequality is proved in \cite{Ledoux}   by using the Kantorovich dual formula  and   Hamilton-Jacobi equations, which are available when $(M,\rr)$ is a   length space as we assumed, see \cite{V2}.
In the following,   we prove estimates \eqref{WP} and \eqref{WP'} by  five steps.

{\bf (a)} For any $r>0$ and $t>0$, consider the regularized empirical measure
\beq\label{MTR} \mu_{t,r}^B:= \mu_t^B \hat P_r=f_{t,r}\mu,\ \ f_{t,r}:=\ff 1 t \int_0^t
\hat p_r(X_s^B,\cdot)\d s.\end{equation}
Note that $\pi:= \ff 1 t\int_0^t \dd_{X_s^B} \times (\dd_{X_s^B} \hat P_r)\d s\in \C(\mu_t^B, \mu_{t,r}^B)$, so that  by the definition of $\W_p$, we obtain
\beg{align*} &\W_p(\mu_{t,r}^B,\mu_t^B)^2 \le \bigg(\int_{M\times M} \rr(x,y)^p\pi(\d x,\d y)\bigg)^{\ff 2 p}  =  \bigg(\ff 1 t \int_0^t \hat P_r\rr(X^B_s,\cdot)^p  (X_s^B)\d s\bigg)^{\ff 2p},\ \ t>0, \ r\in (0,1].\end{align*}
Combining this with  \eqref{A3} and that $\mu$ is an invariant measure of $X_s^B$,    we obtain
\beq\label{00} \E^\mu\big[\W_p(\mu_{t,r}^B,\mu_t^B)^2\big]\le \bigg(\int_{M}  \hat P_r\rr(x,\cdot)^p (x)\mu(\d x)\bigg)^{\ff 2 p} \le k r,\ \ t>0, r\in (0,1].\end{equation}

{\bf (b)} Since $\hat P_s$ has symmetric heat kernel $\hat p_s$,   \eqref{MTR} implies
$ \hat P_s f_{t,r}= \hat P_{\ff{s+r}2} f_{t, \ff{s+r}2}.$
Combining this with  \eqref{Ledoux}, \eqref{A1}, Jensen's inequality,   and that $\mu$ is $\hat P_s$-invariant,    we derive
\beq\label{CG}\beg{split} & \W_p(\mu_{t,r}^B,\mu)^2 \le p^2  \bigg\|\nn  \int_0^\infty  \hat P_{\ff{r+s}2} (f_{t,\ff{r+s}2}-1)\d s\bigg\|_{L^p(\mu)}^2\\
&\le 4^\bb p^2 k^2\bigg(  \int_0^\infty \ff{\e^{-\ll (s+r)/2}}{(s+r)^\bb}  \big\|f_{t,\ff{r+s}2}-1 \|_{L^2(\mu)}  \d s \bigg)^{2}.\end{split}\end{equation}
By   \eqref{MTR}, we obtain
\beg{align*} &\|f_{t,\ff{r+s}2}-1\|_{L^2(\mu)}^2
 =\ff 2{t^2}  \int_0^t\d t_1 \int_{t_1}^t \mu\Big(\big(\hat p_{\ff{r+s}2}(X^B_{t_1},\cdot)-1\big)\big(\hat p_{\ff{r+s}2}(X^B_{t_2},\cdot)-1\big)\Big) \d t_2\\
&= \ff 2 {t^2} \int_0^t\d t_1 \int_{t_1}^t \big[ \hat p_{r+s}(X^B_{t_1},X^B_{t_2})-1\big] \d t_2.\end{align*}
Combining this with \eqref{CG}, we find   a constant $c_1>0$ such that for any $t>0, r\in (0,1]$ and measurable function $h: (0,\infty)\to (0,\infty),$
\beq\label{E1}\beg{split}&   \W_p(\mu_{t,r}^B,\mu)^2 \le   \ff {c_1} {t^2}
  \bigg\{  \int_0^\infty \ff{\e^{-\ff \ll 2(s+r)}}{(s+r)^\bb}  \bigg( \int_0^t\d t_1 \int_{t_1}^t \big[ \hat p_{r+s}(X^B_{t_1},X^B_{t_2})-1\big] \d t_2\bigg)^{\ff 1 2} \d s \bigg\}^{2}\\
  &\le \ff {c_1}{t^2} \bigg(\int_0^\infty \ff{\e^{-\ff \ll 2(s+r)}}{h(s+r)}\d s\bigg)\int_0^\infty \ff{\e^{-\ff \ll 2(s+r)}h(s+r)}{(s+r)^{2\bb}} \d s  \int_0^t\d t_1 \int_{t_1}^t  [ \hat p_{r+s}(X^B_{t_1},X^B_{t_2})-1\big] \d t_2. \end{split}\end{equation}

{\bf (c) } We claim that there exist constants $c',\ll'>0$, such that for all $r>0 $ and  $ t_2>t_1\ge 0,$
\beq\label{E20}  \E^\mu\big[ \hat p_r(X^B_{t_1},X^B_{t_2})-1 \big]\le c' (1\land r)^{-\ff d 4} \Big[r\land 1+1_{\{\aa>0, d'<\infty\}} \{1\land (t_2-t_1)\}^{\ff{d'}{\aa d}}\Big]^{-\ff d 4}\e^{-\ll'(t_2-t_1)}. \end{equation}
Indeed, since   $\mu$ is $P_t^B$-invariant, $\mu$ is also $P_t^{B*}$-invariant, where  $P_t^{B*}$ is the adjoint operator
  of $P_t^B$ in $L^2(\mu)$. By   \eqref{LT} and \eqref{A5}, we obtain that for all $t>0$,
\beq\label{EX*} \ \|P_t^{B*} -\mu \|_{L^2(\mu)}= \|P_t^B -\mu\|_{L^2(\mu)}
 =\big\|\E P_{S_t^B}-\mu\big\|_{L^2(\mu)} \le k \E\big[\e^{-\ll S_t^B}\big]= k\e^{-B(\ll)t}.\end{equation}
Denoting $\ll'=B(\ll),$ noting that
$ \hat p_r(x,\cdot)-1= \int_M \hat p_{\ff r2}(x,y) \big[\hat p_{\ff r 2}(y,\cdot)-1\big]\mu(\d y),$
 by  the Markov property of $X^B_t$,  \eqref{A4} and \eqref{EX*}, we derive
\beq\label{S*} \beg{split}& \E^\mu\big[ \hat p_r(X^B_{t_1},X^B_{t_2})-1 \big]= \int_M P^B_{t_2-t_1} \big[\hat p_r(x,\cdot)-1\big](x) \mu(\d x)\\
&= \int_{M\times M} \hat p_{\ff r 2}(x,y) P^B_{t_2-t_1} \big[\hat p_{\ff r 2}(y,\cdot)-1\big](x)\mu(\d x)\mu(\d y)\\
&\le \bigg(\int_{M\times M} \hat p_{\ff r 2}(x,y)^2\mu(\d x)\mu(\d y)\bigg)^{\ff 1 2} \bigg(\int_M \big\|P^B_{t_2-t_1}[\hat p_{\ff r 2} (y,\cdot)-1]\|_{L^2(\mu)}^2\mu(\d y)\bigg)^{\ff 1 2}\\
&\le k \e^{-\ll'(t_2-t_1)} \int_M \hat p_r(x,x)\mu(\d x) \le k^2 (1\land r)^{-\ff d 2} \e^{-\ll'(t_2-t_1)}.\end{split}\end{equation}
 Next, let $\aa>0$ and $d'<\infty.$  By the Markov property we obtain
\beq\label{E3'} \beg{split} &\E^\mu\big[ \hat p_r(X^B_{t_1},X^B_{t_2})-1 \big]= \int_M \E^x\big[\hat p_r(x, X^B_{t_2-t_1})-1\big] \mu(\d x)\\
&=\int_{M\times M} \big\{\hat p_r(x,y)-1\big\}p^B_{t_2-t_1}(x,y)\mu(\d x)\mu(\d y)\\
&= \int_{M\times M} \hat p_{r}(x,y)  \big\{  p^B_{t_2-t_1}(x,y)-1\big\}  \mu(\d x) \mu(\d y).\end{split}\end{equation}
 Noting that
\beq\label{E21} p^B_{t_2-t_1}(x,y)-1= P_{\ff{t_2-t_1}2}^{B*} \Big\{p^B_{\ff{t_2-t_1}2} (x,\cdot)-1\Big\}(y),\ \ t_2>t_1\ge 0,\end{equation} by  the display after \cite[(3.12)]{WW},  \eqref{A6} implies
\beq\label{A6'} \int_{M\times M}p_t^B(x,y)^2\mu(\d x)\mu(\d y)\le k \E\big[(1\land S_t^B)^{-\ff{d'}2}\big]\le k' (1\land t)^{-\ff{d'}{2\aa}},\ \  t>0\end{equation}
for some constants $k'>0.$
 Combining \eqref{EX*}   with  \eqref{E21}, \eqref{A4} and \eqref{A6},    we find   constants $c_2,c_3>0$ such that
\beg{align*} &\int_{M\times M} \hat p_{r}(x,y)  \big\{  p^B_{t_2-t_1}(x,y)-1\big\}  \mu(\d x) \mu(\d y)\\
 &\le \bigg(\int_{M}\hat p_{2r}(x,x) \mu(\d x) \bigg)^{\ff 1 2} \bigg(\int_M \Big\|P_{\ff{ t_2-t_1}2}^{B*}\Big\{p^B_{\ff{ t_2-t_1}2} (x,\cdot)-1\Big\}\Big\|_{L^2(\mu)}^2\mu(\d x)\bigg)^{\ff 1 2} \\
 &\le c_2 (1\land r)^{-\ff d 4}  \e^{-(t_2-t_1)B(\ll)/2} \bigg(\int_M \Big\| p^B_{\ff{ t_2-t_1}2} (x,\cdot)-1\Big\|_{L^2(\mu)}^2 \mu(\d x)\bigg)^{\ff 1 2} \\
 &\le c_3(1\land r)^{-\ff d 4} \{1\land (t_2-t_1)\}^{-\ff{d'}{4\aa}}   \e^{-(t_2-t_1)\ll' },\ \ \ t_2> t_1\ge 0.\end{align*}
 This together with   \eqref{S*} and \eqref{E3'} implies  \eqref{E20} for some constant $c'>0$.

{\bf  (d) }   When  $ \aa>0$ and $d'<\infty,$ we find   constants $c_4,c_5>0$ such that
 \beg{align*} &\ff 1 t \int_0^t\d t_1\int_{t_1}^t   \e^{-\ll'(t_2-t_1)} \big[(r+s)\land 1  + \{1\land (t_2-t_1)\}^{\ff{d'}{\aa d}}\big]^{-\ff d 4} \d t_2\\
 &\le \ff{c_4}t \int_0^t\d t_1\int_{0}^\infty \e^{-\ll' \theta} \big([(r+s)\land 1]^{\ff {\aa d}{d'} } +  1\land \theta \big)^{-\ff {d'} {4\aa}} \d \theta\le c_5 I(r+s),\end{align*}
 where, since $\ff{d\aa}{d'}(\ff{d'}{4\aa}-1)^+= \ff d 4(1-\ff{4\aa}{d'})^+,$
 \beq\label{I} I(r+s):=     [(1\land (r+s)]^{-\ff d 4(1-\ff{4\aa}{d'})^+}   \Big(1+ 1_{\{d' = 4\aa\}} \log \big[2+ (r+s)^{-1}\big]\Big).\end{equation}
 When $d'=\infty$ or $\aa=0$,  we have
    $I(r+s)=[(1\land (r+s)]^{-\ff d 4}$ and
 $$\ff 1 t \int_0^t\d t_1\int_{t_1}^t   \e^{-\ll' (t_2-t_1)} \big[(r+s)\land 1 \big]^{-\ff d 4} \d t_2\le c_5 I(r+s)$$ holds for some constant $c_5>0.$ So, in any case,
 $$\ff 1 t \int_0^t\d t_1\int_{t_1}^t   \e^{-\ll'(t_2-t_1)} \big[(r+s)\land 1  + 1_{\{\aa>0, d'<\infty\}}\{1\land (t_2-t_1)\}^{\ff{d'}{\aa d}}\big]^{-\ff d 4} \d t_2
 \le c_5 I(r+s).$$
  Combining this  with
  \eqref{E1} and \eqref{E20}, and choosing
  $h(r+s)= \ff{(r+s)^{\bb+\ff d 8}}{\ss{I(s+r)}},$
   we find  a constant   $c_6  >0$ such that
  \beq\label{E*}  \E^\mu[\W_p(\mu_{t,r}^B,\mu)^2]\le \ff{c_6}t \bigg(\int_0^\infty\ff{\e^{-\ll' s}\ss{I(r+s)}}{[1\land (r+s)]^{\bb+\ff d 8}} \d s\bigg)^2.\end{equation}
By \eqref{I} and the definition of $K_{\bb,d,d',\aa},$   we find a constant $c_7>0$ such that
\beg{align*}    \int_0^\infty\ff{\e^{-\ll' s}\ss{I(r+s)}}{[1\land (r+s)]^{\bb+\ff d 8}} \d s
 & \le \int_0^\infty \ff{\e^{-\ll' s}}{[1\land (r+s)]^{K_{\bb,d,d',\aa}}} \Big(1 + 1_{\{d'=4\aa\}} \ss{\log[2+(r+s)^{-1}]}\Big)\d s\\
& \le c_7\eta(r),\ \ \ \ r\in (0,1],\end{align*}
     where
 $$\eta(r):= \beg{cases} 1, &\text{if}\ K_{\bb,d,d',\aa}<1,\\
  \log(2+r^{-1}), &\text{if}\ K_{\bb,d,d',\aa}=1, d'\ne 4\aa,\\
 [\log(2+r^{-1})]^{\ff 3 2}, &\text{if}\ K_{\bb,d,d',\aa}=1, d'= 4\aa,\\
 r^{1-K_{\bb,d,d',\aa}},   &\text{if}\ K_{\bb,d,d',\aa}>1, d'\ne 4\aa,\\
 r^{1-K_{\bb,d,d',\aa}}\ss{\log(2+r^{-1})},   &\text{if}\ K_{\bb,d,d',\aa}>1, d'= 4\aa.\end{cases}$$
 This together with \eqref{E*} and  \eqref{00} implies
\beg{align*} \E^\mu[\W_p(\mu_{t},\mu)^2]\le 2\inf_{r\in (0,1]} \big\{\E^\mu[\W_p(\mu_{t,r}^B,\mu)^2]+ \E^\mu[\W_p(\mu_{t,r}^B,\mu_t^B)^2]\big\}
 \le  c_8 \inf_{r\in (0,1]}  \big\{t^{-1} \eta(r)^2 + r\big\} \end{align*} for some constant $c_8>0$ and all $t>0$. Therefore,
  \eqref{WP}  holds for some constant $c>0$.

{\bf (e) } To prove  \eqref{WP'}, let $t>1$ and $\bar \mu_{t-1}^B:= \ff 1 {t-1} \int_1^t \dd_{X_s^B}\d s,$ so that
$ \mu_t^B= \ff 1 t \int_0^1 \dd_{X_s^B}\d s+ \ff {t-1}t\bar\mu_{t-1}^B.$
Then
$$\W_p(\mu_t^B,\mu)^p\le \ff 1 t\int_0^1  \mu\big(\rr(X_s^B,\cdot)^p  \d s+
\ff{t-1}t \W_p(\bar\mu_{t-1}^B,\mu)^p.$$
By Jensen's inequality and the Markov property, this implies
\beq\label{EN} \beg{split} &\E^x\big[\W_p(\mu_t^B,\mu)^q\big]
 \le \ff{2^{q-1}}{t^q} \E^x \bigg(\int_0^1 \mu\big(\rr(X_s^B,\cdot)^p\big)^{\ff 1 p} \d s\bigg)^q + 2^{q-1} \E^x \big[\W_p(\bar\mu_{t-1}^B,\mu)^q\big]\\
& \le \ff{2^{q-1}}{t^q}  \int_0^1 \E^x\Big[ \mu\big(\rr(X_s^B,\cdot)^p\big)^{\ff q p}\Big] \d s  + 2^{q-1}
\E^{\nu_x}\big[\W_p(\mu_{t-1}^B,\mu)^q\big],\end{split}\end{equation}
where $\nu_x:= p_1^B(x,\cdot)\mu$ is the distribution of $X_1^B$ for $X_0^B=x$. By H\"older's inequality we obtain
\beg{align*} &\E^{\nu_x}\big[\W_p(\mu_{t-1}^B,\mu)^q\big]= \int_M p_1^B(x,y) \E^y \big[\W_p(\mu_{t-1}^B,\mu)^q\big]\mu(\d y)\\
&\le \|p_1^B(x,\cdot)\|_{L^{\ff 2 {2-q}}(\mu)} \bigg(\int_M \E^y\big[\W_p(\mu_{t-1}^B,\mu)^2\big]\mu(\d y)\bigg)^{\ff q 2}
 = \|p_1^B(x,\cdot)\|_{L^{\ff 2 {2-q}}(\mu)} \big(\E^\mu [\W_p(\mu_{t-1}^B,\mu)^2]\big)^{\ff q 2}.\end{align*}
 Combining this with \eqref{EN}, we deduce   \eqref{WP'}.

 \subsection{An extension}

 For some infinite-dimensional models, see for instance \cite{W4}, \eqref{A4} fails for any $d\in (0,\infty)$, but there may be  a decreasing function $\gg: (0,\infty)\to (0,\infty)$ such that
$$ \int_M \hat p_t(x,x)\mu(\d x)\le \gg(t),\ \ t>0.$$
Moreover, in case that $P_t$ is not $L^2$-exponential ergodic, by the weak Poincar\'e inequality which holds for a broad class of ergodic Markov processes,   see \cite{RW01}, we  have
$$\lim_{t\to\infty}  \|P_t-\mu\|_{L^\infty(\mu)\to L^2(\mu)}= 0.$$
To cover these two situations for which Theorem \ref{T} does not apply, we present the following result for the empirical measure $\mu_t$ of the Markov process $X_t$ with semigroup $P_t$.

\beg{thm}\label{TE} Assume   $\eqref{A1}$,  $\eqref{A3}$. If there exist a constant $q\in [1,\infty],q'\in[\ff q {q-1},\infty]$ and a decreasing function $\gg: (0,\infty)\to (0,\infty)$ such that
\beq\label{A5'} \lim_{t\to \infty}  \|P_t-\mu\|_{L^{q'}(\mu)\to L^{\ff{q}{q-1}}(\mu)} =0,\end{equation}
\beq\label{A4'}  \int_M \|\hat p_{\ff r 2}(y,\cdot)\|_{L^q(\mu)}\|\hat p_{\ff r 2}(y,\cdot)\|_{L^{q'}(\mu)} \mu(\d y)\le \gg(r),\ \ \ r>0.\end{equation}
Then there exists a constant $c>0$ such that for any $t>0,$
$$\E^\mu[\W_p(\mu_t,\mu)^2] \le c \inf_{r\in (0,1]} \bigg\{\ff {\int_0^t \|P_{s}-\mu\|_{L^{q'}(\mu)\to L^{\ff q {q-1}}(\mu)}\d s} t \bigg(\int_0^1 \ff{\ss{\gg(r+s)}}{(r+s)^\bb}\d s\bigg)^2 +r\bigg\}.$$\end{thm}

\beg{proof}  Let $B(\ll)=\ll$ so that $X_t^B=X_t, P_t^B=P_t$ and $\mu_t^B=\mu_t$. Noting that $\hat p_r(x,\cdot)= \int_M \hat p_{\ff r 2}(x,y)\hat p_{\ff r 2}(y,\cdot)\mu(\d y),$ by  \eqref{A4'} we obtain
\beg{align*} & \E^\mu[\hat p_r(X_{t_1},X_{t_2})-1] =\int_M (P_{t_2-t_1}-\mu)\hat p_r(x,\cdot)(x)\mu(\d x)\\
&=\int_{M\times M} \hat p_{\ff r 2}(x,y) (P_{t_2-t_1}-\mu) \hat p_{\ff r 2}(y,\cdot)(x)\mu(\d x)\mu(\d y)\\
&\le \int_M \|\hat p_{\ff r 2} (\cdot,y)\|_{L^q(\mu)} \|(P_{t_2-t_1}-\mu) \hat p_{\ff r 2}(y,\cdot)\|_{L^{\ff q {q-1}}(\mu)}\mu(\d y)\\
&\le \|P_{t_2-t_1}-\mu\|_{L^{q'}(\mu)\to L^{\ff q {q-1}}(\mu)} \int_M \|\hat p_{\ff r 2}(y,\cdot)\|_{L^q(\mu)}\|\hat p_{\ff r 2}(y,\cdot)\|_{L^{q'}(\mu)}\mu(\d y) \\
&\le \gg(r) \|P_{t_2-t_1}-\mu\|_{L^{q'}(\mu)\to L^{\ff q {q-1}}(\mu)},\ \ r>0, t_2>t_1.\end{align*}
Combining this with \eqref{A5'}, we find a constant $c_1>0$ such that
$$\int_0^t \d t_1\int_{t_1}^t \E^\mu[\hat p_{r+s}(X_{t_1},X_{t_2}) -1] \d t_2 \le c_1   \gg(r+s) t \int_0^t \|P_{s}-\mu\|_{L^{q'}(\mu)\to L^{\ff q {q-1}}(\mu)}\d s,\ \ r,s>0.$$
So, by \eqref{E1} with
$$h(r+s)= (s+r)^\bb\bigg(\gg(r+s)t \int_0^t \|P_{s}-\mu\|_{L^{q'}(\mu)\to L^{\ff q {q-1}}(\mu)}\d s\bigg)^{-\ff 1 2},$$  we find   constants $c_2,c_3>0$ such that
\beg{align*} &\E^\mu[\W_p(\mu_{t,r},\mu)^2] \le \ff{c_2}t \bigg(\int_0^\infty \ff{\e^{-\ff \ll 2(r+s)} \ss{\gg(r+s)}}{(r+s)^\bb} \d s\bigg)^2\int_0^t \|P_{s}-\mu\|_{L^{q'}(\mu)\to L^{\ff q {q-1}}(\mu)}\d s \\
&\le \ff {c_3}t \bigg(\int_0^1 \ff{  \ss{\gg(r+s)}}{(r+s)^\bb} \d s\bigg)^2\int_0^t \|P_{s}-\mu\|_{L^{q'}(\mu)\to L^{\ff q {q-1}}(\mu)}\d s,\ \ t>0, r\in (0,1].\end{align*}
This together with \eqref{00} and the triangle inequality implies the desired estimate.
 \end{proof}

 To verify Theorem \ref{TE}, we present below a simple example where $P_t$ only has algebraic convergence in $\|\cdot\|_{L^\infty(\mu)\to L^2(\mu)},$ so Theorem \eqref{T} does not apply.

\paragraph{Example 2.3.}  Let $M=[0,1], \rr(x,y)=|x-y|$ and $\mu(\d x)= \d x.$ For any $l\in (2,\infty),$ let $X_t$ be the diffusion process on $M\setminus \{0,1\}$ generated by
$$L:= \big\{x(1-x)\big\}^l\ff{\d^2}{\d x^2} + l \big\{x(1-x)\big\}^{l-1}(1-2x) \ff{\d}{\d x}.$$
Then  Theorem \ref{T} does not apply, but by Theorem \ref{TE} there exists a constant $c>0$ such that for any $t>0$,
\beq\label{CV} \E^\mu[\W_p(\mu_t,\mu)^2]\le c \beg{cases} t^{-1}, &\text{if} \ l\in (2,5), p \in [2, \ff{13-l}4), \\
t^{-1}[\log(2+t)]^3,  &\text{if} \ l\in (2,5], p=\ff{13-l}4, \\
 \big[t^{-1} \log(2+t)]^{\ff 8{4p+l-5}}, &\text{if} \   l\in (2,5], p>\ff{13-l}4,\\
 t^{-\ff 4{l-1}}[\log(2+t)]^2, &\text{if}\ l>5, p=2\\
 t^{-\ff 8 {p(l-1)}},  &\text{if} \ l>5, p>2.\end{cases}\end{equation}

\beg{proof}  We first observe that \eqref{A5} fails, so that Theorem \ref{T} does not apply. Indeed, the Dirichlet form of $L$ satisfies
\beq\label{DF} \EE(f,g)=\int_0^1 \big\{x(1-x)\big\}^l (f'g')(x)\d x,\ \ f,g\in C_b^1(M)\subset \D(\EE).\end{equation}
Let $\rr_L$ be the  intrinsic distance  function to the point $\ff 1 2\in M.$ We find a constant $c_1>0$ such that
$$\rr_L(x)= \bigg|\int_{\ff 1 2}^x \big\{s(1-s)\big\}^{-\ff l 2}\d s\bigg|\ge c_1 \big(x^{1-\ff l 2}+ (1-x)^{1-\ff l 2}\big),\ \ x\in M.$$
Then  for any $\vv>0,$ we have $\mu(\e^{\vv \rr_L})=\infty,$ so that by \cite{AMS},  ${\rm \gap}(L)=0.$
On the other hand, since $L$ is symmetric in $L^2(\mu)$, by \cite[Lemma 2.2]{RW01},    \eqref{A5} implies  the same inequality for $k=1$, so that ${\rm gap}(L)\ge \ll>0$.  Hence,  \eqref{A5} fails.

To apply Theorem \ref{TE}, let $\hat P_t$ be the standard Neumann heat semigroup on $M$ generated by $\DD$.
It is classical that \eqref{A1} and  \eqref{A3}    hold  for
\beq\label{BMM} \bb= \ff 1 2+\ff {p-2}4.\end{equation} Moreover,  there exists a constant $c_2>1$ such that
$$\|\hat P_{\ff r 2}\|_{L^m(\mu)\to L^n(\mu)}\le c_2(1+r^{-\ff {n-m}{2nm}}),\ \ \ 1\le m\le n\le\infty, r>0,$$ so that for $q'= \infty$ and $q>1$,
\beg{align*}  &\|\hat p_{\ff r 2}(y,\cdot)\|_{L^q(\mu)}\|\hat p_{\ff r 2}(y,\cdot)\|_{L^{q'}(\mu)}  \le c_2 \|\hat P_{\ff r 2}\|_{L^1(\mu)\to L^q(\mu)}  (1+ r^{-\ff 1 2}) \\
&\le c_2^2 (1+r^{-\ff {q-1} {2q}}) (1+r^{-\ff 1 2}),\ \ r>0.\end{align*}
Hence, there exists a constant $c_3>0$ such that \eqref{A4'} holds for
$$\gg(r)= c_3  (1+ r^{-\ff{2q-1}{2q}}).$$ Combining this with \eqref{BMM}, we find a constant $k>0$ such that for any $r\in (0,1),$
\beq\label{GMM} \eta(r):= \bigg(\int_0^1 \ff{\ss{\gg(r+s)}}{(r+s)^\bb}\d s\bigg)^2\le k\cdot \beg{cases}
1, &\text{if}\ 1<q<\ff 1 {p-2},\\
[\log(1+r^{-1})]^2, &\text{if}\ 1<q=\ff 1 {p-2},\\
r^{\ff{1-(p-2)q}{2q}}, &\text{if}\ q>  1\lor\ff{1}{p-2}.\end{cases}\end{equation}
Indeed,  \eqref{BMM} implies
$$-\ff{2q-1}{4q}-\bb= -\ff{(p+2)q-1}{4q}\beg{cases}>-1, &\text{if}\  1<q<\ff 1 {p-2},\\
=-1, &\text{if} \ 1<q=\ff 1 {p-2},\\
<-1, &\text{if}\  q>1\lor \ff 1 {p-2},\end{cases}$$
so that we find  constants $k_1,k_2>0$ such that for any $r\in (0,1)$,
$$\int_0^1 \ff{\ss{\gg(r+s)}} {(r+s)^{\bb}}\d s\le k_1 \int_0^1 (r+s)^{-\ff{2q-1}{4q}-\bb}\d s\le k_2\cdot \beg{cases}
1, &\text{if}\ 1<q<\ff 1 {p-2},\\
 \log(1+r^{-1}), &\text{if}\ 1<q=\ff 1 {p-2},\\
r^{\ff{1-(p-2)q}{4q}}, &\text{if}\ q>  1\lor\ff{1}{p-2},\end{cases} $$  which implies  \eqref{BMM}.

To calculate $\|P_t-\mu\|_{L^{q'}(\mu)\to L^{\ff q{q-1}}(\mu)} $  for $q'=\infty$, we apply  the weak Poincar\'e inequality studied in \cite{RW01}. Let
$$M_s=[s,1-s],\ \ \ s\in (0,1/2).$$
Noting that $\mu(\d x)=\d x$ and letting  $\nu(\d x)= \{x(1-x)\}^l\d x,$ we find a constant $c_4>0$ such that
$$\sup_{r\in [s,\ff 12]} \mu([r,1/2])  \nu([s,r]) \le  2^l \sup_{r\in [s,\ff 1 2]} \Big(\ff 1 2 -r\Big)  \big(s^{1-l}- r^{1-l} \big)
\le c_4 s^{1-l},\ \ s\in (0,1/2).$$ By the weighted Hardy inequality \cite{B''},  see for instance \cite[Proposition 1.4.1]{Wbook},
we have
$$\mu(f^2 1_{[s,\ff 1 2]})\le 4 c_4 s^{1-l}\nu(|f'|^2),\ \ f\in C^1([s,1/2]), f(1/2)=0.$$
By symmetry, the same holds for $[\ff 1 2,1-s]$ replacing $[s,\ff 1 2].$ So, according to \cite[Lemma 1.4.3]{Wbook}, see also \cite{Chen}, we derive
$$\mu(f^2 1_{M_s}) \le 4 c_4 s^{1-l} \nu(|f'|^21_{M_s}) +\mu(f1_{M_s})^2,\ \ f\in C^1([s,s-1]).$$
Combining this with \eqref{DF}, for any $f\in C^1_b(M)$ with $\mu(f)=0$, we have $\mu(f1_{M_s})= -\mu(f1_{M_s^c})$ so that
\beg{align*} &\mu(f^2) = \mu(f^21_{M_s^c}) + \mu(f^21_{M_s})\le \mu(f^21_{M_s^c}) +  4c_4 s^{1-l} \EE(f,f)  +\mu(f1_{M_s^c})^2,\\
&\le 4 c_4 s^{1-l}   \EE(f,f) + 2 \|f\|_\infty^2 \mu(M_s^c)^2\le 4 c_4 s^{1-l}\EE(f,f) + 8 s^2 \|f\|_\infty^2,\ \ s\in (0, 1/2).\end{align*}
For any $r\in (0,1),$ let $s= (r/8)^{\ff 1 2}.$ We find a constant $c_5>0$ such that
$$\mu(f^2) \le c_5 r^{-\ff {l-1} 2} \EE(f,f) +r\|f\|_\infty^2,\ \ r\in (0,1), \mu(f)=0, f\in C^1_b(M).$$
By \cite[Corollary 2.4(2)]{RW01}, this implies
$$\|P_t-\mu\|_{L^\infty(\mu)\to L^2(\mu)}=\|P_t-\mu\|_{L^2(\mu)\to L^1(\mu)}\le c_5 (1+t)^{-\ff 2 {l-1}},\ \ t>0$$ for some constant $c_5>0.$
Since $ P_t $ is contractive in $L^n(\mu)$ for any $n\ge 1$,  this together with the interpolation theorem  implies
$$\|P_t-\mu\|_{L^\infty(\mu)\to L^{\ff q {q-1}}(\mu)}\le c_6 (1+t)^{-\ff {4(q-1)}{q(l-1)}},\ \ t>0.$$
Noting that $q'=\infty,$ we find a constant $k>0$ such that
\beq\label{GGT} \GG(t):= \ff 1 t \int_0^t \|P_s-\mu\|_{L^{q'}(\mu)\to L^{\ff q {q-1}}(\mu)}\d s\le k\beg{cases}
 t^{-1},  &\text{if}\  l\in (2,5), q>\ff 4{5-l},   \\
t^{-1} \log(2+t),  &\text{if}\ l=5, q=\infty,\\
(1+t)^{-\ff{4 }{l-1}},  &\text{if}\ l>5, q=\infty.\end{cases}\end{equation}
We now prove the desired estimates case by case.

(1) Let $l\in (2,5)$ and $p\in [2,\ff{13-l}4).$ Taking  $q\in (\ff 4{5-l}, \ff 1 {p-2})$ in \eqref{GMM} and \eqref{GGT},  we obtain
$$\inf_{r\in (0,1]} \big\{\eta(r) \GG(t)+r\big\}\le k \inf_{r\in (0,1]}\big\{t^{-1}+r\big\}= kt^{-1}.$$
So, the desired estimate follows from Theorem \ref{TE}.

(2) Let $l\in (2,5]$ and $p=\ff{13-l}4.$ Taking $q= \ff 4{5-l}=\ff 1 {p-2}$ in  \eqref{GMM} and \eqref{GGT} we find a constant $c>0$ such that
$$\inf_{r\in (0,1]} \big\{\eta(r) \GG(t)+r\big\}\le k \inf_{r\in (0,1]}\big\{t^{-1}[\log(2+t)][\log (1+r^{-1})]^2+r\big\}\le c t^{-1}[\log (2+t)]^3.$$
This implies the desired estimate  according to  Theorem \ref{TE}.

(3) Let $l\in (2,5]$ and $p>\ff{13-l}4.$ We have  $q:= \ff 4 {5-1}>\ff 1 {p-2},$ so that   \eqref{GMM} and \eqref{GGT} imply
$$\inf_{r\in (0,1]} \big\{\eta(r) \GG(t)+r\big\}\le k \inf_{r\in (0,1]}\big\{t^{-1}[\log(2+t)] r^{-\ff{4p+l-13} 8}+r\big\}\le c \big[t^{-1}\log (2+t)\big]^{-\ff{4p+l-5}8} $$
for some constant $c>0$, which implies the desired estimate by Theorem \ref{TE}.

(4) Let $l>5$ and $p=2$. By taking $q=\infty$ in  \eqref{GMM} and \eqref{GGT}, we find a constant $c>0$ such that
$$\inf_{r\in (0,1]} \big\{\eta(r) \GG(t)+r\big\}\le k \inf_{r\in (0,1]}\big\{t^{-\ff 4{l-1}}[ \log(1+r^{-1})]^2 +r\big\}\le c  t^{-\ff 4{l-1}}[\log(2+t)]^2.$$
By Theorem \ref{TE},  the desired estimate holds.

(5) Let $l>5$ and $p>2$. By taking $q=\infty$ we find a constant $c>0$ such that  \eqref{GMM} and \eqref{GGT} imply
$$\inf_{r\in (0,1]} \big\{\eta(r) \GG(t)+r\big\}\le k \inf_{r\in (0,1]}\big\{t^{-\ff 4{l-1}} r^{-\ff{p-2}2} +r\big\}\le c  t^{-\ff 8{p(l-1)}}$$
for some consatnt $c>0.$ Hence
   the desired estimate holds according to Theorem \ref{TE}.

\end{proof}

\section{Subordinate stochastic Hamiltonian systems  }

Consider the following degenerate SDE for $X_t=(X_t^{(1)},X_t^{(2)})$ on $\R^{n+m}=\R^{n}\times\R^{m}$ ($n,m\ge 1$ may be different):
\beq\label{E2} \beg{cases} \d X_t^{(1)}= \kk Q X_t^{(2)}\d t,\\
\d X_t^{(2)}= \ss 2 \,\d W_t -\big\{Q^*(\nn V) (X_t^{(1)})+\kk X_t^{(2)}\big\}\d t,\end{cases}\end{equation}
where $W_t$ is the $m$-dimensional Brownian motion,  $Q\in\R^{n\otimes m},$ $\kk>0$ is a constant, and $V\in C^2(\R^{n})$ satisfies
\beq\label{EP}\sup_{x_1\in\R^{n}} \ff{ \|\nn^2 V(x_1)\| }{ 1+|\nn V(x_1)|}<\infty,\ \   \int_{\R^{n}}    |\nn V(x_1)|^2 \e^{-V(x_1)}\d x_1<\infty.\end{equation}
Let
$$\mu_V(\d x_1):=\ff{\e^{-V(x_1)}\d x_1}{\int_{\R^{n}}\e^{-V(x_1)}\d x_1},\ \ \scr N_\kk(\d x_2):=\Big(\ff{\kk}{2\pi}\Big)^{\ff {m} 2}\e^{-\ff \kk  2 |x_2|^2}\d x_2.$$
Then the SDE \eqref{E2} is well-posed, and the solution has invariant probability measure
\beq\label{MU} \mu(\d x_1,\d x_2):=     \mu_V(\d x_1)\scr N_\kk(\d x_2).\end{equation}
Recall that for a metric space $(M,\rr)$,
$$B(x,r):=  \big\{y\in M:\ \rr(x,y)  <r\big\},\ \ \ x\in M, r>0.$$
 We will verify  $(A_1)$ and \eqref{A6} for the present model by using the following assumption.

\emph{\beg{enumerate} \item[$(A_2)$] $QQ^*$ is invertible, $\eqref{EP}$ holds, and there exist   constants $k>0$ and $n'\ge n$ such that
\beq\label{H1}   \nn^2 V\ge -k I_{n},\ \  \int_{\R^{n}} \ff {\mu_V(\d x_1)} {\mu_V(B(x_1,r))} \le k r^{-n'},\  \ r\in (0,1],\end{equation}
\beq\label{H2} \mu_V(f^2)\le k \mu_V(|\nn f|^2),\ \ f\in C_b^1(\R^{n}),\ \mu_V(f)=0.\end{equation}
\end{enumerate}}
 We have the following result.

\beg{thm}\label{T2} Assume $(A_2)$,   let $B\in {\bf B}^\aa$ for some $\aa\in [0,1],$  and let $\rr(x,y)=|x-y|$ for $x,y\in \R^{n+m}.$ Let $\xi_t(K)$ be in $\eqref{KK}$ for
\beq\label{KK*} K:= \beg{cases} \ff 1 2 +\ff{n'+2m}4 -\ff{\aa(n'+2m)}{2(3n'+2m)}, &\text{if}\ \|\nn^2V\|_\infty<\infty\\
\ff 1 2 +\ff{n'+2m}4, &\text{if}\ \|\nn^2V\|_\infty=\infty.\end{cases} \end{equation}
Then  there exists a constant $c>0$ such that
\beq\label{K1}   \E^\mu \big[\W_2(\mu_t^B,\mu)^2\big]\le c  \xi_t(K),\ \ \ t>0.\end{equation}
If $\|\nn^2 V\|_\infty<\infty,$ then  for any $t\ge 2$ and $x\in \R^{n+m},$
\beq\label{K2} \big[\E^x  \W_2(\mu_t^B,\mu)\big]^2 \le c \xi_t(K)\E^x \bigg[\int_0^1   |X_s^B|^2 \d s + \ff 1 {\mu\big(B(x_1, (1\land S_1^B)^{\ff 3 2})\times B(x_2,(1\land S_1^B)^{\ff 1 2})\big)}     \bigg].\end{equation}
\end{thm}

To prove this result, we first present a dimension-free Harnack inequality for   the following more general model:
\beq\label{E0'} \beg{cases} \d X_t^{(1)}= \big\{A X_t^{(1)} +QX_t^{(2)}\big\}\d t,\\
\d X_t^{(2)}= Z_t(X_t)\d t+\si_t \d W_t,\end{cases}\end{equation}
where $Q$ and $W_t$ are   in \eqref{E2},   $A \in\R^{n\otimes n},$ and
$$\si: [0,\infty)\to \R^{m\otimes m},\ \ \ Z: [0,\infty)\times \R^{n+m}\to\R^{m}$$ are measurable  such that the following conditions hold:
\emph{\beg{enumerate}\item[$(A_3)$] There exist  a constant $k>0$ and an integer $k_0\ge 0$ such that
  $$\sup_{t\ge 0} \|\si_t^{-1}\|_\infty+\sup_{t\ge 0, x\ne y} \ff{|Z_t(x)-Z_t(y)|} {|x-y|}\le k,\ \ \
  {\rm Rank}[A^i Q: 0\le i\le k_0]= n. $$ \end{enumerate} }

\beg{lem}\label{LNN} Assume $(A_3)$, and let $P_t$ be the Markov semigroup associated with $\eqref{E0'}.$  Then for any $p\in (1,\infty)$, there exists a constant $c(p)>0$ such that
$$|P_t f(y)|^p\le \big(P_t |f|^p(x) \big)\exp\bigg[ \ff{c|x_1-y_1|^2}{(1\land t)^{4k_0+3}}+ \ff{c|x_2-y_2|^2}{(1\land t)^{4k_0+1}} \bigg]$$ holds for all $t>0$ and $ x=(x_1,x_2),y=(y_1,y_2)\in \R^{n }\times \R^m.$
\end{lem}

\beg{proof}  By Jensen's inequality, we only need to prove for $t\in (0,1].$ The proof is refined from that of \cite[Lemma 3.2]{W17}. Let $X_t$ solve \eqref{E0'} for $X_0=x$, and for fixed $t_0\in (0,1],$ let $Y_t$ solve the following SDE
with $Y_0=y$:
$$ \beg{cases} \d Y_t^{(1)}= \big\{A Y_t^{(1)} +QY_t^{(2)}\big\}\d t,\\
\d Y_t^{(2)}= \Big\{Z_t(X_t)+\dfrac{x_2-y_2}{t_0} +\dfrac{\d }{\d t}\big[t(t_0-t)Q^*\e^{(t_0-t)A^*}b_{t_0}\big]\Big\}\d t+\si_t \d W_t,\ \ t\in [0,t_0],\end{cases}$$
where
\beg{align*} & b_{t_0}:= Q_{t_0}^{-1} \bigg\{\e^{t_0A} (x_1-y_1)+ \int_0^{t_0} \ff{t_0-s}{t_0}\e^{(t_0-s)A}Q^* (x_2-y_2)\d s\bigg\},\\
&Q_{t}:= \int_0^{t} s(t-s)\e^{(t-s)A} QQ^*\e^{(t-s)A^*}\d s,\ \ t>0.\end{align*}
By \cite[(3.2) and (3.3)]{W17}, we have $X_{t_0}=Y_{t_0}$ and
 \beq\label{*P2} \sup_{t\in [0,t_0]} |X_t-Y_t|\le c  |x-y|,\ \  x,y\in \R^{n+m}\end{equation} holds
 for some constant $c>0.$

 According to the proof of \cite[Theorem 4.2]{WZ}, the rank condition in $(A_3)$  implies
 $$\|Q_{t_0}^{-1}\|\le c_1 t_0^{-2k_0-3},\ \ \ t_0\in (0,1]$$  for some constant $c_1>0.$  Then there exists a constant $c_2>0$ such that
 \beq\label{*P1} |b_{t_0}|\le c_2 t_0^{-2k_0-3}\big(|x_1-y_1|+t_0 |x_2-y_2|\big).\end{equation}
Combining this with the first condition in $(A_3)$, we see that
 $$\psi_t:= \si_t^{-1} \bigg(Z_t(X_t)-Z_t(Y_t) + \ff{x_2-y_2}{t_0} +\ff{\d}{\d t} \Big\{t(t_0-t) Q^*\e^{(t_0-t)A^*}b_{t_0}\Big\}\bigg) $$
 satisfies
 \beq\label{*P3} \sup_{t\in [0,t_0]} |\psi_t|_\infty^2 \le c_4\Big(\ff{|x_1-y_1|^2}{t_0^{4k_0+4}}+\ff{|x_2-y_2|^2}{t_0^{4k_0+2}}\Big),\ \ t_0\in (0,1]\end{equation}
 for some constant $c_4>0.$
So, by Girsanov's theorem, under the probability measure $R\d\P,$ where
 $$R:= \exp\bigg[-\int_0^{t_0} \<\psi_t,\d W_t\>-\ff 1 2 \int_0^{t_0} |\psi_t|^2\d t\bigg],$$ the process  $(Y_t)_{t\in [0,t_0]}$ is a weak solution to   \eqref{E0'} with initial value $y$.
Combining this together with   $X_{t_0}=Y_{t_0}$ as observed above, we find a constant $c_5>0$ depending on $p$  such that
 \beg{align*}& |P_{t_0}f(y)|^p =|\E[R f(Y_{t_0})]|^p= |\E[R f(X_{t_0})]|^p\le \big(\E[R^{\ff p{p-1}}]\big)^{p-1} \E[|f|^p(X_{t_0})]\\
 &\le \e^{c_5 \int_0^{t_0} \|\psi_t\|_\infty^2 \d t}P_{t_0}|f|^p(x).\end{align*}
 By \eqref{*P3}, this implies the desired Harnack inequality.
\end{proof}

\beg{proof}[Proof of Theorem \ref{T2}] Let $p=2, M:=\R^{n+m}.$  To apply Theorem \ref{T}, let $\hat X_t$ be the  diffusion process generated by
\beq\label{HL} \hat L:=\DD-(\nn H)\cdot\nn,\end{equation}
where
\beq\label{H} H(x): = V(x_1)+ \ff { \kk}2 | x_2|^2,\ \ x=(x_1,x_2)\in \R^{n+m}.\end{equation}
    In the following,  we    verify  $(A_1)$  and Definition \ref{DEF}    for
 \beq\label{GG} \bb=\ff 1 2 ,\ \  d= n'+2 m,\ \ d'=  \beg{cases} 3n'+2m, &\text{if}\ \|\nn^2 V\|_\infty<\infty,\\
 \infty, &\text{otherwise}.\end{cases}\end{equation}

(a)  Verify \eqref{A1}. By \eqref{MU},  \eqref{H2} and the Poincar\'e inequality for the Gaussian measure $\scr N_\kk$, we find a constant $C>0$ such that
\beq\label{PP2} \mu(f^2)\le C \mu(|\nn f|^2) +\mu(f)^2,\ \ f\in C_b^1(\R^{n+m}).\end{equation}
 Consequently, ${\rm gap}(\hat L)\ge C^{-1}>0$, so that  Proposition \ref{P1} implies \eqref{A1}  for $p=2$ and $\bb=\ff 1 2.$

 (b) Verify \eqref{A3} and \eqref{A5}.   By \eqref{H1} and \eqref{H}, there exists a constant $c_1>0$ such that
\beq\label{HHS}  \nn^2 H   \ge -c_1 I_{n+m}.\end{equation}
 Then as in \eqref{BO},  we find a constant $c_2>0$ such that
 $$\hat L|x-\cdot|^2\le c_2 (1+ |\nn V(x_1)|^2+|x-\cdot|^2),\ \ x\in \R^{n+m}.$$ This implies
 \beq\label{B'} \hat P_t|x-\cdot|^2(x)=\E^x|x-\hat X_t|^2 \le c_2\big(1+\nn V(x_1)|^2\big)t\e^{c_2t},\ \ x\in\R^{n+m}, \ t>0.\end{equation}
  Combining this with \eqref{EP} and \eqref{MU}, we verify condition \eqref{A3} for $p=2$ and some constant $k>0.$
 Moreover,  according to  \cite{GS},   \eqref{EP} and \eqref{H2} imply  \eqref{A5}   for some constants $k,\ll>0.$

 (c)  Verify \eqref{A4}.
 According to   \cite[Theorem 2.4.4]{Wbook}, by \eqref{HHS} we find a constant $c_3>0$ such that
 \beq\label{AA}\hat p_r(x,x)\le \ff {c_3} {\mu(B(x,\ss r))},\ \ x\in \R^{n+m}, r\in (0,1].\end{equation}
 Combining this with \eqref{MU}, \eqref{H1}  and   $B(x,\ss r) \subset B(x_1, \ss {r/2})\times B(x_2, \ss{r/2})$, we find a constant $c_4>0$ such that
\beq\label{A0} \int_{\R^{n+m}} \ff {\mu(\d x)} {\mu(B(x,\ss r))} \le  c_4   r^{-\ff {n'}2} \int_{\R^{m}} \ff {\e^{-\kk |x_2|^2} \d x_2 } {\int_{B(x_2,\ss{r/2})}\e^{-\kk |y_2|^2}\d y_2}, \ \  r\in (0,1].\end{equation}
By the same argument leading to \eqref{CC},    we find a   constant $c_5>0$ such that
 \beq\label{**P}  \int_{\R^{m}} \ff{\e^{-\kk |x_2|^2}\d x_2}{\int_{B(x_2,r)}\e^{-\kk |y_2|^2}\d y_2} \le     c_5r^{-2m},\ \ r\in (0,1]. \end{equation}
Combining this with \eqref{A0}, we  find a constant $c >0$ such that
\beq\label{A0'} \int_{\R^{n+m} }\ff{\mu(\d x)}{\mu(B(x,\ss r))}\le  c  r^{-\ff{n'+2m} 2},\ \ \ r\in (0,1].\end{equation}
Since  $\hat p_t(x,x)$ is decreasing in $t>0$,  this together with \eqref{AA} implies  \eqref{A4} for $d= n'+2m.$

(d)    To estimate $d'$, we assume $\|\nn^2 V\|_\infty<\infty.$  By  Lemma \ref{LNN}   for $p=2$, where $k_0=0$ holds for the present model,    we find a constant $c_6>0$ such that
 $$|P_tf(x)|^2\le (P_tf^2(y))\e^{c_6 (1\land t)^{-3}|x_1-y_1|^2+ c_6(1\land t)^{-1}|x_2-y_2|^2},\ \ t>0, x,y\in \R^{n+m}.$$
  Choosing $f:= p_t(x,\cdot)\land l,$ we derive
\beg{align*} &\bigg(\int_{\R^{n+m}}(p_t(x,\cdot)\land l)^2\d \mu \bigg)^2\e^{-c_6 (1\land t)^{-3}|x_1-y_1|^2- c_6(1\land t)^{-1}|x_2-y_2|^2} \\
&\le P_t (p_t(x,\cdot)\land l)^2(y),\ \ l\ge 1.\end{align*}
 Integrating both sides with respect to $\mu(\d y)$ and noting that $\mu$ is $P_t$-invariant, we obtain
\beg{align*} &\int_{\R^{n+m}}(p_t(x,\cdot)\land l)^2\d \mu \le \ff 1 {\int_{\R^{n+m}} \e^{-c_6 (1\land t)^{-3}|x_1-y_1|^2- c_6(1\land t)^{-1}|x_2-y_2|^2}\mu(\d y)}\\
&\le \ff {\e^{2c_6}}{\mu\big(B(x_1,(1\land t)^{\ff 3 2})\times B(x_2, (1\land t)^{-\ff 1 2})\big)}.\end{align*}
 Letting $l\to\infty$ we arrive at
 \beq\label{HTE} \int_{\R^{n+m}}p_t(x,\cdot)^2\d \mu \le   \ff {\e^{2c_6}}{\mu\big(B(x_1,(1\land t)^{\ff 3 2})\times B(x_2, (1\land t)^{-\ff 1 2})\big)},\ \ t>0, x\in\R^{n+m}.\end{equation}
 This together with  \eqref{H1}  and \eqref{**P}  yields
\beg{align*} &\int_{\R^{n+m}\times\R^{n+m}}p_t(x,y)^2\mu(\d x)\mu(\d y)\\
&\le  \e^{2c_6}  \int_{\R^{n+m}} \ff{\mu(\d x)} {\mu\big(B(x_1,(1\land t)^{\ff 3 2})\times B(x_2, (1\land t)^{-\ff 1 2})\big)}\le c (1\land t)^{-\ff{3n'+2m}2}\end{align*}
 for some constant $c>0$.
 Therefore,   \eqref{A6} holds for $d'=  3n'+2m.$

(e) For $K$ in \eqref{KK*}, $\bb, d,d'$ in \eqref{GG} and $\aa\in [0,1]$, we have
$$d' \ge  3n'+2m\ge 5  >4\aa,\ \ \
  K_{\bb,d,d',\aa} =K.$$
Then \eqref{K1}  follows from \eqref{WP}.

Next,  by \eqref{PP2} we have $\mu(\e^{c|\cdot|})<\infty$ for some constant $c>0,$ see for instance \cite{AMS}.
Combining this with  \eqref{B}, we find a constant $c_7>0$ such that
$$\E^x[\mu(\rr(X_s^B,\cdot)^2)]= \E^x[\mu(|X_s^B-\cdot|^2)]\le c_7(1+\E^x |X_s^B|^2).$$
Moreover, \eqref{HTE} implies
$$ \|p_1^B(x,\cdot)\|_{L^2(\mu)}^2 \le \E\bigg[\ff {\e^{2c_6}}{\mu\big(B(x_1,(1\land S_1^B)^{\ff 3 2})\times B(x_2, (1\land S_1^B)^{-\ff 1 2}\big)}\bigg].$$
Then \eqref{K2}  follows from   \eqref{WP'}
 \end{proof}

\paragraph{Example 3.1.}  Consider \eqref{E2} with invertible    $QQ^*$, and let $V$ be in Example 2.1.
Then  for any $\aa\in [0,1]$ and $B\in {\bf B}^\aa,$    there exists a constant $c>0$ such that for any $t\ge 1$,
\beq\label{K3}   \E^\mu [\W_2(\mu_t^B,\mu)^2]\le c \beg{cases} t^{-\ff{(2\tau-1)(3\tau n +  2\tau m-m )}{(\tau n+ 2\tau m-m)][3\tau n +(m-\aa) (2\tau-1)] }}, &\text{if}\ \ff 1 2< \tau\le 1,\\
t^{-\ff {2\tau-1} {\tau n+m(2\tau-1)}}, &\text{if}\ \tau>1.\end{cases} \end{equation}
 When $\aa>0$ and $\tau\in (\ff 1 2,1],$  there exists a constant $\vv>0$ such that for any $  t\ge 1$ and $x\in \R^{n+m},$
\beq\label{K4} \big[\E^x  \W_2(\mu_t^B,\mu)\big]^2 \le c  t^{-\ff{(2\tau-1)(3\tau n +  2\tau m-m )}{(\tau n+ 2\tau m-m)][3\tau n +(m-\aa) (2\tau-1)] }}  \e^{\kk   |x_2|^2 +  (1+\theta|x_1|^2)^\tau-\vv  |x_2|- \vv |x_1|^\tau}.
\end{equation}

\beg{proof}    (1) As explained in the proof of Example 2.1 that
$(A_2)$ holds  for $n'=\ff{2\tau n}{2\tau-1}$ and some constant $k>0$.  So, $K$ defined in \eqref{KK*} satisfies $K>1$.  It is easy to see that  $\|\nn^2 V\|_\infty=\infty$ for $\tau>1$ while$\|\nn^2 V\|_\infty<\infty$ for $\tau\in (\ff 1 2,1]. $ Then
estimate \eqref{K3}   follows from Theorem \ref{T2}.

(2)  Let $\aa>0$ and $\tau\le 1.$ We find a constant $c_1>0$ such that
\beq\label{NM} \sup_{t\in [0,1] } \E^x[|X_t|^2]  \le c_1(1+|x|^2).\end{equation}
Next,  similarly to \eqref{SO}, there exists a constant $c_2>0$ such that
$$\mu\big(B(x_1,r^{\ff 3 2})\times B(x_2,r^{\ff 1 2})\big)\ge c_2 r^{\ff {3n+m} 2 } \e^{c_2 r^{\ff 3 2} |x_2|+ c_2r^{\ff \tau 2}|x_1|^\tau-\kk |x_2|^2-(1+\theta |x_1|^2)^\tau},\ \ r\in (0,1], x\in\R^{n+m}. $$
Combining this with \eqref{HTE}, we find a constant $c_3 >0$ such that
\beq\label{*8}\beg{split} &\int_{\R^{n+m}} p_t(x,\cdot)^2\d\mu \le    \ff {\e^{2c_6}}{\mu\big(B(x_1,(1\land t)^{\ff 3 2})\times B(x_2, (1\land t)^{-\ff 1 2}\big)}\\
&\le c_3^2  \ff{\exp[ \kk   |x_2|^2 +  (1+\theta|x_1|^2)^\tau-c_2 r^{\ff 3 2} |x_2|- c_2r^{\ff \tau 2}|x_1|^\tau]}{(1\land t)^{\ff{3n+m}2}},\ \ \ t>0, x\in \R^{n+m}.\end{split}\end{equation}
Consequently,
\beg{align*} &\E^x[|X_t|^2]=\int_{\R^{n+m}} |y|^2 p_t(x,y)\mu(\d y)\le  \bigg(\mu(|\cdot|^4) \int_{\R^{n+m}}p_t(x,y)^2\mu(\d y)\bigg)^{\ff 1 2}\\
&\le c_3 \ff{\exp[\ff \kk 2 |x_2|^2 +\ff 1 2 (1+\theta|x_1|^2)^\tau-\ff 1 2 c_2 r^{\ff 3 2} |x_2|-\ff 1 2  c_2r^{\ff \tau 2}|x_1|^\tau]}{(1\land t)^{\ff{3n+m}4}},\ \ t>0, x\in \R^{n+m}.\end{align*}
This together with \eqref{NM} yields
\beq\label{*9}\sup_{s\in [0,1]} \E^x[|X_s^B|^2]\le \sup_{t\ge 0} \E^x[|X_t|^2]  \le c \e^{\ff \kk 2 |x_2|^2 +\ff 1 2 (1+\theta|x_1|^2)^\tau},\ \ x\in \R^{n+m}.\end{equation}
Moreover, when $\aa>0$,  the second inequality in \eqref{A6'} implies
$$\E\big[(1\land S_1^B)^{-\ff {3n+m } 2}\big]<\infty, $$ which together with  \eqref{*8} yields
$$\E^x\bigg[ \ff {1}{\mu\big(B(x_1,(1\land S_1^B)^{\ff 3 2})\times B(x_2, (1\land S_1^B)^{-\ff 1 2}\big)}\bigg]\le c\e^{ \kk   |x_2|^2 +  (1+\theta|x_1|^2)^\tau-\vv  |x_2|- \vv |x_1|^\aa}$$
for some constant $\vv>0.$
Combining this   with   \eqref{*9},
we deduce   the \eqref{K4}   from that in Theorem \ref{T2}.

\end{proof}

\section{Subordinate spherical velocity Langevin diffusions}

In this section, we consider the following degenerate SDE on $M:=\R^n\times \mathbb S^{n-1} (n\ge 2)$:
\beq\label{ES} \beg{cases}\d X_t^{(1)} = X_t^{(2)}\d t,\\
\d X_t^{(2)}= - \ff 1 {n-1} \big(I_n - X_t^{(2)}\otimes X_t^{(2)}\big) \nn V(X_t^{(1)}) \d t +\si \big(I_n-X_t^{(2)}\otimes X_t^{(2)}\big) \circ \d W_t,\end{cases}\end{equation}
where $V\in C^2(\R^n),$  $W_t$ is an $n$-dimensional Brownian motion,   $\si>0$ is a constant, and  $\circ \d  $ is the Stratonovich differential.
The solution of \eqref{ES} is called  the  spherical velocity Langevin diffusion process generated by
$$L:= \ff { \si^2}2 \DD^{(2)}  + x_2 \cdot\nn^{(1)} - (\nn^{(2)}\Phi) \cdot \nn^{(2)},$$
where $\DD^{(2)}$ and $\nn^{(2)}$ are the Laplacian and gradient on $\mathbb S^{n-1}$ respectively,   $\nn^{(1)}$ is the gradient on $\R^n$, and
$$\Phi(x):= \ff 1 {n-1} \big(\nn^{(1)} V(x_1)\big)\cdot x_2,\ \ \ x=(x_1,x_2)\in M. $$
See \cite{GS} and references therein for the background of this model.

Let $\rr$ be the  Riemannian distance on $M:= \R^{n}\times \mathbb S^{n-1},$   let $V$ satisfy $(A_2)$, and  let
$$\mu(\d x):= \mu_V(\d x_1) \LL(\d x_2),$$
where $\LL$ is the normalized volume measure on $\mathbb S^{n-1}.$ We have  the following result.

\beg{thm}\label{T3} Let $V$ satisfy  $(A_2)$ and   let $B\in {\bf B}. $
Then there exists a constant $c>0$ such that for any $t\ge 2,$
$$  \E^\mu\big[ \W_2(\mu_t^B,\mu)^2\big]\le c  t^{-\ff 2{n'+n-1}}, $$
$$\big[\E^x  \W_2(\mu_t^B,\mu)\big]^2 \le c     t^{-\ff 2{n'+n-1}}\E^x \bigg[\int_0^1   |X_s^B|^2 \d s + \ff 1 {\mu(B(x, (1\land S_1^B)^{\ff 3 2}))} \bigg],\ \ x\in \R^{n}\times \mathbb S^{n-1}.$$
In particular, for $V$ given in Example $2.1$,  these estimates hold for $n'=2n.$
 \end{thm}

\beg{proof}   According to \cite[Theorem 1.1]{GS}, $(A_2)$ implies \eqref{A5}   for some constants $k,\ll>0.$
 Let $\DD$ be the Laplacian on $M:=\R^n\times \mathbb S^{n-1}.$ To apply Theorem \ref{T}, we choose the reference symmetric diffusion process generated by
$$\hat L:= \DD -\{\nn^{(1)}V(x_1)\}\cdot\nn^{(1)},$$
which is symmetric in $L^2(\mu)$.  As shown in the proof of Example 3.1,  \eqref{H1} holds for   $n'=2n.$   So,  by Theorem \ref{T} for $d'=\infty$, it suffices to  verify
  \eqref{A1}, \eqref{A3} and \eqref{A4}  for  $p=2,\    \ \bb=\ff 1 2 ,\  d= n'+n-1.$

   By $(A_2)$ and   the compactness  of  $\mathbb S^{n-1}$, the Bakry-Emery curvature of $\hat L$ is bounded  below by a constant, and there exists a constant $\ll>0$ such that
$$\mu(f^2)\le \ff 1 {2\ll} \mu(|\nn f|^2) +\mu(f)^2,\ \ f\in C_b^1(M).$$
Then as explained in steps (a)-(b) in the proof of Theorem \ref{T2},     \eqref{A1} and \eqref{A3} hold  for $p=2,\bb=\ff 1 2$ and some constants $k>0.$

  By \cite[Theorem 2.4.4]{Wbook} and that  the  Bakry-Emery curvature of $\hat L$ is bounded from below, we find a constant $c_1>0$ such that
 $$\hat p_r(x,x)\le \ff {c_1} {\mu(B(x,\ss r))},\ \ x\in M, r\in (0,1].$$
 Noting that
 $$\LL(B(x_2,r))\ge c_2 r^{n-1},\ \ \ r\in (0,1]$$ holds for some constant $c_2>0,$ as explained in step (c) in the proof of Theorem \ref{T2},
 we derive   \eqref{A4} for $d= n'+n-1.$

 \end{proof}

\section{Subordinate Wright-Fisher type diffusion  processes}

For $1$-dimensional Wright-Fisher diffusions, the convergence rate has been derived for the empirical measure with respect to the Wasserstein distance induced by
the intrinsic distance. However, for  higher dimensional Wright-Fisher type diffusion  processes, the intrinsic distance is less explicit, so that it is hard to apply the framework introduced in
\cite{W23}. Below, we consider higher dimensional Wright-Fisher type diffusion  processes using Wasserstein distance induced by the Euclidean distance rather the intrinsic distance, so that
the framework introduced in the present paper works well.

 Let $2\le N\in \mathbb N,$ consider
 $$\DD^{(N)}:= \Big\{x\in [0,1]^N: |x|_1:=\sum_{i=1}^N x_1\le 1\Big\}.$$
Given $q=(q_i)_{1\le i\le N+1}\in [1,\infty)^{N+1}$,
the  Dirichlet distribution $\mu$ with parameter $q$  is a probability measure on $\DD^{(n)}$ defined as follows:
\beq\label{WF} \beg{split} &\mu(\d x):= 1_{\DD^{(N)}}(x) h(x)\d x,\ \ h(x):= \ff{\GG(|q|_1)}{\prod_{i=1}^{N +1} \GG(q_i)}  \prod_{i=1}^{N+1} x_i^{q_i-1} ,\\
&  |q|_1:=\sum_{i=1}^{N+1}q_i,  \ \ x_{N+1}:=1-|x|_1, \ \ x\in \DD^{(N)}.\end{split}\end{equation}
 The Dirichlet distribution   arises naturally in Bayesian inference as conjugate priors for categorical distribution, and also arises in population genetics describing the distribution of  allelic frequencies, see for instance
 \cite{CM} and references within.

 Let  $X_t$ be either the Wright-Fisher diffusion with mutation    generated by
  $$\sum_{i,j=1}^N \big(x_i\dd_{ij}-x_ix_j\big) \pp_{x_i}\pp_{x_j} + \sum_{i=1}^N (q_i-|q|_1x_i) \pp_{x_i},$$
 or the   diffusion process generated by
 $$\sum_{i=1}^N \Big\{x_i(1-|x|_1)\pp_{x_i}^2 +\big(q_i(1-|x|_1)-q_{N+1}x_i\big)\pp_{x_i}\Big\}.$$Then the associated Dirichlet form is determined by
\beq\label{DR} \EE(f,g):=\beg{cases} \mu\big(\sum_{i,j=1}^N (x_i\dd_{ij}-x_ix_j)(\pp_{x_i}f)(\pp_{x_j}g)\big), &\text{the \ first\ case,}\\
 \mu\big(\sum_{i,j=1}^N  x_i(1-|x|_1)(\pp_{x_i}f)(\pp_{x_j}g)\big), &\text{the \ second\ case,}\end{cases}
 \end{equation}
 for any $f,g\in C_b^1(\DD^{(N)}).$
In both cases,  the  Poincar\'e inequality
 \beq\label{PI} \mu(f^2)\le \ff 1 {\ll } \EE(f,f),\ \ f\in C_b^1(\DD^{(N)}), \mu(f)=0\end{equation}
 holds for some constant $\ll >0$, see
   \cite{ST} and \cite{FMW} for the value of the    largest constant $\ll $ (i.e. the spectral gap).

 \beg{thm}\label{TNN}   Let $p\in [2,\infty)$,  $B\in {\bf B}^\aa$ for some $\aa\in [0,1],$  and $\rr(x,y)=|x-y|, x,y\in \DD^{(N)}.$ Let
 $\xi_t(K)$ be in $\eqref{KK}$  for
 \beg{align*} K:= \ff 1 2 +\ff d 4 -\ff{\aa d}{2 d'},\ \ d:=   |q|_1-1,\ \ d':= 4\sum_{i=1}^N q_i +2q_{N+1}-2.\end{align*}
 Then there exists a constant $c>0$ such that
  $$  \sup_{x\in \DD^{(N)}} \E^x[\W_p(\mu_t^B,\mu)^2]\le c \xi_t(K),\ \ t>0.  $$
 \end{thm}

To apply Theorem \ref{T}, let    $\hat P_t$ be the Neumann semigroup on $\DD^{(N)}$ generated by
 $$\hat L:=\DD+ (\nn\log h)\cdot\nn,$$
 where $h$ is in \eqref{WF}.
 The associated Dirichlet form is determined by
 $$\hat \EE(f,f)= \mu(|\nn f|^2),\ \ \ f\in C_b^1(\DD^{(N)}).$$ To verify
 \eqref{A1}, we first present the following lemma on the super Poincar\'e inequality of $\hat\EE$.

 \beg{lem}\label{LS}   For the above $\hat \EE$ and $\mu$, there exists a constant $c>0$ such that
 \beq\label{SPI} \mu(f^2)\le r \hat\EE(f,f)+c \big(1+r^{-\ff{|q|_1-1}2}\big)\mu(|f|)^2,\ \ f\in \D(\hat\EE), r>0.\end{equation}
 \end{lem}
 \beg{proof} (a) We  follow the idea of \cite{WZ19}. For any $s>1$ and $r>0,$  let
\beg{align*} &D_s:= \big\{x\in \DD^{(N)}:\ \phi(x):=(1-|x|_1)^{-1}\le s\big\},\\
&\ll(s):= \inf\big\{\hat\EE(f,f):\ f\in C^1(\DD^{(N)}), f|_{D_s}=0, \mu(f^2)=1\big\},\\
&s_r:= \inf\big\{s>1: \ll(s)\ge 8 r^{-1}\big\},\end{align*}
so that  $\phi(x):= (1-|x|_1)^{-1} $ satisfies
$$ h(s):= \sup_{D_s}\hat\EE(\phi,\phi)= N s^4,\ \ s>1.$$
According to \cite[Theorem 2.1]{WZ19}, if
\beq\label{SPP} \mu(f^2)\le r \hat \EE(f,f)+\bb_s(r) \mu(|f|)^2,\ \ r>0, f\in C^1(\DD^{(N)}), f|_{D_s^c}=0\end{equation} holds for some
$\bb_s: (0,\infty)\to (0,\infty),$ then
 there exists a constant $c_1>0$ such that
\beq\label{SPI'} \mu(f^2)\le r \hat\EE(f,f)+\bb(r)\mu(|f|)^2,\ \ f\in \D(\hat\EE), r>0 \end{equation} holds for
\beq\label{BBT} \bb(r):= c_1 + \big(2+8Nrs_r^2\big)\bb_{3s_r}\Big(\ff r{8+4Nrs_r^2}\Big),\ \ r>0.\end{equation}

(b) Let $\EE$ be in \eqref{DR}. By \cite[Lemma 3.2]{WZ19}, there exist constants $c_2,s_0>0$ such that
$$  \inf\big\{\EE(f,f):\ f\in C^1(\DD^{(N)}), f|_{D_s}=0, \mu(f^2)=1\big\}\ge c_2 s,\ \ s\ge s_0.$$
Noting that $\hat\EE(f,f)\ge s\EE(f,f)$ holds for $f\in C^1(\DD^{(N)}), f|_{D_s}=0,$ this implies
$$\ll(s)\ge c_2 s^2,\ \ s\ge s_0.$$ Hence, we find a constant $c_3>0$ such that
\beq\label{SR} s_r\ge c_3 (1+r^{-\ff 12}),\ \ \ r>0.\end{equation}

(c) To estimate $\bb_s(r)$ in \eqref{SPP}, we first consider the following product  probability measure $\tt\mu$ on $[0,1]^N$:
$$\tt\mu(\d x):=\prod_{i=1}^N \mu_i(\d x_i),\ \ \mu_i(\d x_i):= q_ix_i^{q_i-1}  \d x_i,\ \ 1\le i\le N.$$
For any $r\in (0,\ff 1 2]$ and    $I:= [a,b]\subset [0,1]$   with $\mu_i(I)= r$,  we intend to prove
\beq\label{R1}\mu_i^\pp((\pp I) \setminus \{0,1\})=\mu_i^\pp(\{a,b\}\setminus \{0,1\})\ge q_i r^{1-q_i^{-1}},\end{equation}
where $\mu_i^\pp(\{s\}):= q_i s^{q_i-1}$ for $s\in [0,1]$ is the boundary measure induced by $\mu_i.$
We have
$$b^{q_i} = \mu([0,b]) \ge \mu(I)= r,$$ so that $b\ge r^{q_i}.$
If $b<1$, then
$$\mu_i^\pp((\pp I )\setminus \{0,1\})\ge \mu_i^\pp(\{b\})= q_i b^{q_i-1}\ge q_ir^{1-q_i^{-1}}.$$
When $b=1$, we have
$$1-a^{q_i} =\mu_i([a,1]) =\mu_i(I)=r,$$
which implies $a\ge (1-r)^{q_i^{-1}}$ and hence, for $r\in (0,\ff 1 2]$,
$$\mu_i^\pp((\pp I )\setminus \{0,1\})\ge \mu_i^\pp(\{a\})= q_i a^{q_i-1}\ge q_i(1-r)^{1-q_i^{-1}}\ge q_ir^{1-q_i^{-1}}.$$
In conclusion, \eqref{R1} holds for any $r\in (0,\ff 1 2]$ and interval $I\subset [0,1]$ with $\mu_i(I)=r$, so that
$$\kk(r):=\inf_{\mu_i(I)\le r}\ff{\mu_i^\pp((\pp I)\setminus \{0,1\})}{\mu_i(I)}\ge q_i r^{-q_i^{-1}},\ \ r\in (0,1/2].$$
This implies
$$\kk^{-1}(2r^{-\ff 1 2}):=\sup\big\{r'\in (0,1/2]:\ \kk(r')\ge 2 r^{-\ff 1 2}\big\}\ge \ff 12\land \big\{q_i^{q_i}r^{\ff {q_i}2}\big\},\ \ r>0.$$
According to \cite[Theorem 3.4.16(1)]{W05}, we find a constant $c_4>0$ such that this implies
$$\mu_i(f^2)\le r \mu_i(|f'|^2)+ c_4(1+r^{-\ff{q_i}2})\mu_i(|f|)^2,\ \ f\in C^1([0,1]), r>0, 1\le i\le N.$$
By   \cite[Proposition 2.2]{WZ19}, we derive
\beq\label{SPTT} \tt\mu (f^2)\le r \tt\mu (|\nn f|^2)+ c_5\Big(1+r^{-\ff 1 2 \sum_{i=1}^Nq_i}\Big)\tt\mu (|f|)^2,\ \ f\in C^1([0,1]^N), r>0.\end{equation}
Now, given $f\in C^1(\DD^{(N)})$ with $f|_{D_s^c}=0,$ take
$$g(x):= f(x)(1-|x|_1)^{\ff{q_{N+1}-1}2},\ \ x\in \DD^{(N)}.$$ By \eqref{SPTT}, we find constants $c_6,c_7>0$ such that
\beg{align*} \mu(f^2)&=c_6 \tt\mu(g^2) \le c_6 r_1 \tt\mu(|\nn g|^2) + c_6c_5\Big(1+r_1^{-\ff 1 2 \sum_{i=1}^Nq_i}\Big)\tt\mu(|g|)^2\\
& \le c_7 r_1 \mu(|\nn f|^2) + c_7 s^2 r_1\mu(f^2  ) + c_7 \Big(1+r_1^{-\ff 1 2 \sum_{i=1}^Nq_i}\Big) s^{q_{N+1}-1} \mu(|f|)^2,\ \ r_1>0.\end{align*}
For any $r>0,$ by choosing
$$r_1=\ff 1 {2c_7} \big(r \land s^{-2}\big) $$
 in the above inequality, we find a constant $c_8>0$  such that  \eqref{SPP} holds for
 $$\bb_s(r):= c_8 \Big(1+ (r\land s^{-2})^{-\ff 1 2 \sum_{i=1}^Nq_i}\Big)s^{q_{N+1}-1},\ \ s>1, r>0.$$
Combining this with  \eqref{SR}, we find a constant $c>0$ such that $\bb(r)$ in \eqref{BBT} satisfies
$$ \bb(r)\le c (1+r^{-\ff {|q|_1-1} 2}),\ \ r>0,$$
and hence \eqref{SPI} follows from \eqref{SPI'}.
 \end{proof}

 \beg{proof}[Proof of Theorem \ref{TNN}]
 Note that $\rr$ is bounded. Moreover, by \eqref{N1} below, $p_1^B=\E[p_{S_1^B}]$, and the second inequality in \eqref{A6'}, we see that
  $p_1^B$ is also  bounded.  Then    the desired estimate follows from \eqref{WP'} for $q=2$. Noting that  $q_i\ge 1$ implies
 $$d'= 4\sum_{i=1}^{N}q_i+2q_{N+1}-2\ge 4\sum_{i=1}^N q_i\ge  8>4\aa,$$
    by Theorem \ref{T},  it suffices to verify $(A_1)$ and \eqref{A6}  for    the given constants $d, d'$ and
 $$\bb:=  \ff 1 2 + \ff{d(p-2)}{4p},$$   for which we have $K_{\bb,d,d',\aa}=K.$

  (1)  Let $h$ be in \eqref{WP}. Since $q_i\ge 1,$ we have
    $ \nn^2 \log h\le 0,$  so that the Bakry-Emery curvature of $\hat L$ is nonnegative. Since $\DD^{(N)}$ is convex,
    according to \cite[Theorem 3.3.2(11)]{Wbook}, we have
     \beq\label{GR} |\nn \hat P_t f|\le \ff 1 { \ss t} (\hat P_t |f|^2)^{\ff 1 2},\ \ t>0.\end{equation}
Since
   the Dirichlet form $\hat \EE$ of $\hat P_t$ is larger than   $\EE$,    \eqref{PI}  implies the same inequality for $\hat \EE$, so that
  \beq\label{PCC} \|\hat P_t-\mu\|_{L^2(\mu)}\lor \|P_t-\mu\|_{L^2(\mu)}\le \e^{-\ll t},\ \ t\ge 0.\end{equation} Then \eqref{A5} holds.
    Moreover, it is easy to see that \eqref{SPI} implies the Nash inequality
    $$\mu(f^2)\le c_1\big(1+\hat\EE(f,f)\big)^{\ff d{d+2}},\ \ f\in \D(\hat\EE), \mu(|f|)=1$$
    for some constants $c_1>0$, so that by  \cite[Corollary 2.4.7]{Davies},
   \beg{align*}&\|\hat P_t\|_{L^1(\mu)\to L^\infty(\mu)}\le \|\hat P_{\ff t 2}\|_{L^1(\mu)\to L^2(\mu)}\|\hat P_{\ff t2}\|_{L^2(\mu)\to L^\infty(\mu)}\\
   &= \|\hat P_{\ff t 2}\|_{L^1(\mu)\to L^2(\mu)}^2\le c (1\land t)^{-\ff d 2},
    \ \ t>0\end{align*}
     holds for some constant $c>0$. Combining this with the interpolation theorem,  we find  a constant $c_2>0$ such that
   \beq\label{LD} \|\hat P_t\|_{L^{q_1}(\mu) \to L^{q_2}(\mu)}\le c_2 \Big(1+t^{-\ff {d(q_2-q_1)}{2q_1q_2}}\Big),\ \ t>0, 1\le q_1\le q_2\le \infty.\end{equation}
      Taking $q_1=1$ and $q_2=\infty$ we derived \eqref{A4}, while  choosing $q_1=2,q=p$ we deduce \eqref{A1} with   $\bb=\ff 1 2+\ff{d(p-2)}{4p}$
  from  \eqref{GR} and \eqref{PCC}.

  On the other hand,   by \cite[Proof of (3.1)]{WZ19},   the super Poinar\'e inequality
   $$\mu(f^2)\le r\EE(f,f) + c_3(1+r^{1-2\sum_{i=1}^N q_i -q_{N+1}})\mu(|f|)^2,\ \ r>0, f\in \D(\EE)$$
   holds for some constant $c_3>0$. Noting that $2\sum_{i=1}^N q_i +q_{N+1}-1=\ff{d'}2,$   as explained above, we find a constant $c>0$ such that
   \beq\label{N1}\|P_t\|_{L^1(\mu)\to L^\infty(\mu)}\le c (1\land t)^{-\ff{d'}2},\ \ t>0,\end{equation} so that \eqref{A6} holds.

 (2) It remains to verify \eqref{A3}.
 Since $q_i-1\ge 0,$ for any $x,y\in \DD^{(N)}$, we have
 \beg{align*}& \sum_{i=1}^N (y_i-x_i) \Big(\ff {q_i-1}{x_i}- \ff{q_{N+1}-1}{1-|y|_1}\Big) = \sum_{i=1}^N \ff{(q_i-1)(y_i-x_i)}{y_i} +\ff{(q_{N+1}-1) (|x|_1-|y|_1)}{1-|y|_1}\\
 &\le \sum_{i=1}^{N+1}(q_i-1)<\infty.\end{align*}
 So, for any $2\le n\in \mathbb N,$ we find   constants $c_1(n),c_2(n)>0$ such that
 \beg{align*}\hat L |x-\cdot|^n(y)&\le c_1(n) |x-y|^{n-2} +    n|x-y|^{n-2} \sum_{i=1}^N (y_i-x_i) \Big(\ff {q_i-1}{x_i}-\ff{q_{N+1}-1}{1-|y|_1}\Big)\\
 &\le c_2(n)|x-y|^{n-2},\ \ \ x,y\in \DD^{(N)}.\end{align*}
This implies
$$\E^x|x-\hat X_t^x|^n\le c_2(n) \int_0^t \E^x|x-\hat X_s^x|^{n-2}\d s,\ \ t\ge 0.$$
By induction in $n\ge 2$, we find constants $\{c(n)>0\}_{n\ge 2}$ such that
$$ \E^x|x-\hat X_t^x|^n\le c(n) t^{\ff n 2},\ \ t\ge 0, x\in\DD^{(N)}, n\ge 2.$$
 In particular, \eqref{A3} holds for some constant $k>0.$

 \end{proof}

\section{ Subordinate   stable like processes  }

In this section, we consider the subordinate stable like  process in a   connected  closed  smooth  domain $M\subset \R^n, n\in\mathbb N.$

Let $0\le V\in C^2(M)$ such that
$$\mu=\mu_V:= \ff{\e^{-V(x)}\d x }{\int_{M}\e^{-V(x)}\d x} $$ is a probability measure on $M$. Let $\dd' \in (0,2)$, and let $X_t$ be the Markov process on $M$ associated with the $\aa'$-stable like
Dirichlet form
\beq\label{BF} \beg{split} \EE(f,g):= &\,\int_{M\times M} \ff{(f(x)-f(y))(g(x)-g(y))}{|x-y|^{n+\aa'}}\mu(\d x)\mu(\d y),\\
&\qquad   f,g\in  \D(\EE):=\{f\in L^2(\mu):\EE(f,f)<\infty\}. \end{split}\end{equation}
Let $\mu_t^B$   be the empirical measure of $X_t^B:=X_{S_t^B},$ where $B\in {\bf B}.$

\subsection{Bounded $M$}

When $M$ is bounded, we have the following result for the above defined $\mu_t^B$.

 \beg{thm}  \label{TN}  Assume that $M$ is bounded. Let $\rr(x,y)=|x-y|$ for $x,y\in M, $  let $B\in {\bf B}^\aa$ for some $\aa\in (0,1],$  and  let $\xi_t$ be in $\eqref{xi}$ for
\beq\label{PP1} \bb:= \ff 1 2 +\ff{n(p-2)}{4p},\ \ d:=n,\ \ d':=\ff {2n} {\aa'}.\end{equation}
  Then for any  $p\in [2,\infty),$ there exists a constant $c>0$ such that
 $$\sup_{x\in M}\E^x [\W_p(\mu_t^B,\mu)^2]\le c \xi_t,\ \ \ t\ge 2.  $$
 \end{thm}

\beg{proof} Let $\hat L=\DD+\nn V$ be   equipped with the Neumann boundary condition.   It is well-known that \eqref{A1}-\eqref{A4} hold for  some constant $\ll>0$ and constants $\bb$ and $d'$ in \eqref{PP1}.
  According to \cite[Theorem 1.1]{CK}, there exists a constant $k>0$ such that
\beq\label{ENN} \sup_{x,y\in M} p_t(x,y)\le k(1\land t)^{-\ff n \aa'},\ \ t>0.\end{equation}
Although only $V=0$ is considered  in \cite{CK}, this estimate also holds for $V\in C^2(M)$. Indeed,     according to \cite[Corollary 2.4.7]{Davies},  \eqref{ENN} is equivalent to the Nash inequality
$$\mu(f^2)\le c \big(1+\hat\EE(f,f)^{\ff n{n+\aa'}}\big),\ \  f\in\D(\hat\EE), \mu(|f|)=1$$
for some constant $c>0$. If this inequality   holds for $V=0$, then it also holds for bounded $V$, which is the case when $V\in C^2(M)$ for compact  $M$.

By \eqref{ENN},  \eqref{A6} holds for $d'= \ff {2n} \aa'$, and the generator $L$ of $P_t$ has discrete spectrum. Since   the Dirichlet form $\EE$  is irreducible, the discreteness of the spectrum implies the existence of a spectral gap,
so that   \eqref{A5}   holds for  some constant $\ll>0.$
 Since  $\rr$  is  bounded, and  by the second inequality in \eqref{A6'},   \eqref{ENN} implies
 $$\|p_1^B\|_\infty\le k \E[(1\land S_1^B)^{-\ff n \aa'}]=:c_0<\infty.$$
 Combing this with $\|\rr\|_\infty:=\sup_{x,y\in M}|x-y|<\infty$ since $M$ is bounded, and applying \eqref{WP'} for $q=2$, we obtain
 $$\E^x[\W_p(\mu_t^B,\mu)^2] \le \ff 2 {t^2} \|\rr\|_\infty^2 + 2 c_0 c\xi_{t -1},$$
 which implies the desired estimate,  since   there exists a constant $c'>0$ such that   $t^{-1}\lor \xi_{t-1}\le c' \xi_t$ holds for $ t\ge 2.$   \end{proof}

\paragraph{Remark 6.1.}  To show the optimality of Theorem \ref{TN},  we simply consider
  $p=2,\aa=1, \aa'=2$,  and denote $\mu_t^B=\mu_t.$ Then Theorem \ref{TN} implies
$$\sup_{x\in M} \E^x[\W_2(\mu_t,\mu)^2] \le c  \beg{cases} t^{-1}, &\text{if}\ n\le 3,\\
t^{-1} \{\log (2+t)\}^3, &\text{if}\ n=4,\\
t^{-\ff 2 {n-2}},    &\text{if}\ n\ge 5,\end{cases} $$
which is sharp   except for $n=  4  $ where the exact   order is $t^{-1}\log (2+t)$ according to    \cite{WZ23}. This minor loss is caused by the fact that in general $P_t\ne \hat P_t$, so that in  proof of Theorem \ref{T}
we can not combining them together by using  the semigroup property as in \cite{WZ23}.

\subsection{$M=\R^n.$ }
Assume that $V$ satisfies the following conditions  for  some constants $k>0,d\ge n$:
\beg{align}  &\nn^2 V\ge -kI_n,  \ \ \e^{-V}\in C_b^2(\R^n),\label{HV1} \\
  &\liminf_{|x|\to\infty} \ff 1 {|x|} \<\nn V(x),x\> >0,\label{HV2}\\
 &\int_{\R^n} \ff{\mu(\d x)}{\mu(B(x,r))}\le k r^{-d},\ \ \ r\in (0,1],\label{HV3}\\
 &\inf_{x\in\R^n} \ff{\e^{V(x)}}{(1+|x|)^{n+\aa'}}>0,  \label{HV4}\\
&\lim_{r\to\infty} \inf_{|x|\ge r-1}r^{n+\aa'-1}  \big\{|\nn \e^{-V}(x)|  +r^{-1} \e^{-V(x)}\big\}=0.\label{HV5}
\end{align}
We have the following result for $\mu_t^B$ defined for the subordinate process of the jump process with Dirichlet form \eqref{DF}.

\beg{thm}\label{T6} Assume  $\eqref{HV1}$-$\eqref{HV5},$ and let   $B\in {\bf B}.$ Then the following assertions hold.
\beg{enumerate} \item[$(1)$] There exists a constant $c>0$ such that
\beq\label{LS1} \E^\mu [\W_2(\mu_t^B,\mu)^2]\le c t^{-\ff 2 d},\ \ \ t\ge 1.\end{equation}
\item[$(2)$] Let $B\in {\bf B}^\aa$ for some $\aa\in (0,1]$, and there exist constants $k, d'>1$ such that
\beq\label{LS} \|P_t\|_{L^1(\mu)\to L^\infty(\mu)} \le k (1\land t)^{-\ff{d'} 2},\ \ \ t\in (0,1],\end{equation}  then there exists a constant $c>0$ such that
\beq\label{LS2} \sup_{x\in \R^n} \E^x[\W_2(\mu_t^B,\mu)^2] \le c \xi_t,\ \ t\ge 2,\end{equation}
where $\xi_t$ is in $\eqref{xi}$ for $\bb=\ff 1 2,$ $d$ in $\eqref{HV3}$ and $d'$ in $\eqref{LS}$. \end{enumerate}
\end{thm}
\beg{proof} (1) Let $\hat L=\DD-\nn V$. Then \eqref{HV2} implies
$$\limsup_{|x|\to\infty}\hat L |\cdot|(x) <0,$$
so that \cite[Corollary 1.4]{W99} implies  $\gap(\hat L)>0, $ so that by Proposition \ref{P1},  \eqref{A1} holds for $p=2$ and $\bb=\ff 1 2.$
Moreover, \eqref{A3} and \eqref{A4} have been verified in steps (b) and (c) in the proof of Theorem \ref{T2}, respectively. Finally, by \cite[Theorem 1.1(1)]{WW15}, \eqref{HV4} and \eqref{HV5} imply the Poincar\'e inequality for $\EE$, so that \eqref{A5} holds for $k=1$ and some constant $\ll>0$. Therefore, the estimate \eqref{LS1} follows from \eqref{WP} with $\bb=\ff 1 2,\aa=1$ and $d'=\infty$.

(2) When \eqref{LS} holds, by the second inequality in \eqref{A6'}, we have
$$\|p_1^B\|_\infty\le k \E[(1\land S_1^B)^{-\ff{d'}2}]<\infty.$$
 Moroever, \eqref{HV2} implies
$ \mu(|\cdot|^l)<\infty$ for all $l\in (1,\infty)$. Combining this with  \eqref{LS}, for large enough $l$ we find constants $c_1,c_2>0$ such that
\beg{align*} &\sup_{x\in \R^n}\int_0^1 \E^x|X_s^B|^2\d s \le \E \int_0^1 \big\|P_{S_s^B}|\cdot|^2\big\|_\infty\d s\\
&\le  \mu(|\cdot|^{2l})^{\ff 1 l} \E \int_0^1 \|P_{S_s^B}\|_{L^l(\mu)\to L^\infty(\mu)} \d s \\
&\le c_1 \int_0^1 \E\big[(1\land S_s^B)^{-\ff{d'}{2l}}\big]\d s \le c_2 \int_0^1 s^{-\ff{d'}{2l\aa}}\d s<\infty.\end{align*}
Therefore, \eqref{LS2} follows from \eqref{WP'} with $q=2$.

\end{proof}

\paragraph{Example 6.1.} Let $0\le V\in C^2(\R^n)$ be in Example 2.1 for some $\tau>\ff 1 2$. Let
$B\in {\bf B}^\aa$ for some $\aa\in (0,1],$ and let $\xi_t(K)$ be in \eqref{KK} for
$$K_\dd :=\ff 1 2 + \ff{\tau n(\dd n +\dd\aa'-\aa\aa')}{2\dd (2\tau-1)(n+\aa')},\ \ \dd>2.$$ Then for any $\dd>2$ there exists a constant $c>0$ such that
\beq\label{LSS}   \sup_{x\in \R^n} \E^x[\W_2(\mu_t^B,\mu)^2]\le c \xi_t(K_\dd),\ \ t\ge 1.\end{equation}

\beg{proof}
Conditions \eqref{HV1}, \eqref{HV2},   \eqref{HV4}, and \eqref{HV5} trivially hold. Next, as shown in the proof of \eqref{CC},
 \eqref{HV3} holds for $d=\ff{2\tau n}{2\tau-1}$. Moreover, by \cite[Corollary 1.5]{WW15}, for any $\dd>2$ there exists a constant $k>0$ such that
$$ \|P_t\|_{L^1(\mu)\to L^\infty(\mu)} \le k (1\land t)^{-\ff{\dd(n+\aa')}{\aa'} },\ \ \ t\in (0,1],$$
so that \eqref{LS} holds for $d'=\ff{2\dd(n+\aa')}{\aa'},$ which satisfies $d'>2\aa$. It is easy to see that for $\bb=\ff 1 2$ and $d,d'$ given above we have
$K_\dd=K_{\bb, d,d',\aa},$ so that
   \eqref{LSS} follows from \eqref{LS2}.

\end{proof}
 \paragraph{Acknowledgement.}  The author would like to thank Dr. Jie-Xiang Zhu for corrections, and the referees for helpful comments.
\paragraph{Data availability statements.} The manuscript has no associated data.
\paragraph{Declarations.} The work is supported in
 part by   the National Key R\&D Program of China (No. 2022YFA1006000, 2020YFA0712900) and NNSFC (11921001).
The authors have no relevant financial or non-financial interests to disclose.

\end{document}